\documentclass[11pt,twoside]{amsart}

\usepackage{amsfonts}
\usepackage{amssymb}

\usepackage{amscd}
\usepackage{amsthm}
\usepackage[dvips]{graphicx}

\usepackage[all]{xy}

\newtheorem{Theo1}{{Theorem}}
\newtheorem{Lemma1}{{Lemma}}[section]
\newtheorem{Def1}{{Definition}}
\newtheorem{Prop1}[Lemma1]{{Proposition}}
\newtheorem{Claim1}[Lemma1]{{Claim}}
\newtheorem{Rem1}[Lemma1]{{Remark}}
\newtheorem{Cor1}[Lemma1]{{Corollary}}
\newtheorem{Ex1}[Lemma1]{{Example}}

\newenvironment{Lemma}{\begin{Lemma1}}{\end{Lemma1}}

\newenvironment{Prop}{\begin{Prop1}}{\end{Prop1}}
\newenvironment{Rem}{\begin{Rem1}\rm}{\end{Rem1}}
\newenvironment{Theorem}{\begin{Theo1}}{\end{Theo1}}

\newenvironment{Cor}{\begin{Cor1}}{\end{Cor1}}
\newenvironment{Claim}{\begin{Claim1}}{\end{Claim1}}

\newcommand{\dar}{\downarrow}
\newcommand{\lra}{\longrightarrow}
\newcommand{\lla}{\longleftarrow}
\newcommand{\ra}{\rightarrow}
\newcommand{\epi}{\lra \kern-.8em\ra}
\newcommand{\sdp}{\times\kern-.2em\vrule height1.1ex depth-.05ex}

\newcommand{\Q}{{\mathbb Q}}
\newcommand{\N}{{\mathbb N}}

\newcommand{\F}{{\mathbb F}}

\newcommand{\ol}{\overline}
\newcommand{\Sn}{{\mathfrak S}}
\newcommand{\cal}{\mathcal}
\newcommand{\Z}{{\mathbb Z}}
\newcommand{\dickebox}{{\vrule height5pt width5pt depth0pt}}

\normalbaselines
\title{The Ring of Polynomial Functors of Prime Degree}

\author{Alexander Zimmermann}
\address{Universit\'e de Picardie,
D\'epartement de Math\'ematiques et CNRS UMR 7352, 33 rue St Leu,
F-80039 Amiens Cedex 1, France}
\email{alexander.zimmermann@u-picardie.fr}
\thanks{This research was supported by a grant ''PAI alliance'' from
the Minist\`ere des Affaires \'Etrang\`eres de France and
the British Council. The author acknowledges
support from STIC Asie of the Minist\`ere des Affaires \'Etrang\`eres de France}
\date{April 15, 2013; revised July 26, 2013}
\subjclass[2010]{16H10; 20C30; 20J06; 55R40}
\keywords{Polynomial functors; Green orders; Brauer tree algebras;
Schur algebras; Recollement diagram; Representation type}

\begin{document}

\begin{abstract} Let $\hat\Z_p$ be the ring of $p$-adic integers.
We prove in the present paper that the category of polynomial functors
from finitely generated free abelian groups to $\hat \Z_p$-modules
of degree at most $p$ is equivalent to
the category of  modules over a particularly
well understood ring, called Green order. That this is the case
was conjectured by Yuri Drozd.
\end{abstract}

\maketitle

\section*{Introduction}

Polynomial functors attained a lot of interest in recent years by
at least two major discoveries. First, in \cite{HLS} Henn, Lannes and
Schwartz showed that the category of analytic functors from the category
of finite dimensional vector spaces to to the category of vector
spaces over the same field of characteristic $p$ is equivalent to
the category of unstable modules over the mod $p$ Steenrod algebra
modulo nilpotent objects. Second, Franjou, Friedlander, Scorichenko and
Suslin
in \cite{FFSS}, Friedlander and Suslin in \cite{friedlander},
Touz\'e \cite{Touze} as well as Touz\'e and
van der Kallen \cite{TouzevanderKallen}
use strict polynomial functor to prove the finite generation
of cohomology of group schemes and to compute $Ext$-groups
of modules over general linear groups.
More recently Djament and Vespa studied stable
homology of orthogonal, symplectic and unitary groups using some category of
polynomial functors \cite{Vespa0, Vespa1, Vespa2, Djament, DjamentVespa, DjamentVespa2}.
For definitions and more ample remarks of these concepts
we refer to section \ref{generalities}.

The category ${\mathcal A}_R$ of polynomial functors $\Z-free\lra
R-mod$ for a commutative ring $R$ is a classical object in
algebraic topology (cf Eilenberg, MacLane \cite{eilenberg}).
Let ${\mathcal A}_R^n$ be the full subcategory of
at most degree $n$ polynomial functors in ${\mathcal A}_R$.
Quadratic functors were characterized by Baues \cite{Baues}
as modules over a particular algebra. Baues, Dreckmann,
Franjou and Pirashvili show in \cite{BDFP} that
${\mathcal A}_{R}^n$ is a module category of  finitely generated
$R$-algebra $\Gamma_R^n$ as well.
This description was used by Drozd to show
in \cite{drozdquadratic} that ${\mathcal A}_{\hat\Z_2}^2$ and
in \cite{drozdcubic} that ${\mathcal A}_{\hat\Z_3}^3$ are two very
explicitly given classical orders over
$\hat\Z_2$ and $\hat\Z_3$ respectively, whose representation theory
is completely understood. In particular each of them admit only a
finite number of indecomposable lattices. Here, and in the sequel, we
denote by $\hat\Z_p$ the ring of $p$-adic integers, and by
$\F_q$ the field with $q$ elements.
These orders were introduced by Roggenkamp
in \cite{roggreen}. Recall
that an $R$-order over an integral domain $R$ is an
$R$-algebra $\Lambda$, finitely generated projective as $R$-modules and
so that $K\otimes_R\Lambda$ is a semisimple $K$-algebra, for $K$
being the field of fractions of $R$.

Drozd
conjectures at the end of \cite{drozdcubic} that
${\mathcal A}_{\hat\Z_p}^p$ should be equivalent to the module category
of a particular Green order $\Lambda_p$ over $\hat\Z_p$
for all primes $p$.
Drozd proves the case $p=2$ and $p=3$ by explicitly
associating the
generators of the ring given by Baues or Baues, Dreckmann,
Franjou and Pirashvili respectively to matrices in
the corresponding matrix rings and computes the
kernel and the image of the so-defined mapping. The relations in
\cite{BDFP} are sufficiently involved so that going beyond
$p=3$ by this method seems to be not realistic.

\medskip

In this paper we prove Drozd's conjecture.
Our method is conceptual. We develop a recollement diagram
of ${\mathcal A}_{\F_p}^n$ by ${\mathcal A}_{\F_p}^{n-1}$
and the module category of the group ring
$\F_p\Sn_n$, analogous to the one described by
Schwartz~\cite[Section 5.5]{unstable} for functors
$\F_p-mod\lra\F_p-mod$. This recollement diagram for
${\mathcal A}_{\F_p}^n$ may be of independent interest
since it is completely general. It actually appears already in
work of Pirashvili~\cite{Pirashvilirussian},
as it was indicated by the referee.
Comparison of these two diagrams gives many informations.
A second ingredient then is the study of various $Ext$-groups
between simple functors, using work of Franjou, Friedlander,
Scorichenko and Suslin \cite{FFSS}. The third main ingredient is
the explicit projective functor mapping to the reduction
modulo $p$ functor. It should be noted that we do not actually
use the ring defined by Baues, Dreckmann, Franjou
and Pirashvili in \cite{BDFP}. We just use that there is an algebra
which is finitely generated, so that the Krull-Schmidt property,
lifting of idempotents and similar properties are valid for
${\mathcal A}_{\hat\Z_p}^n$. For this reason we do not give a
Morita bimodule between the Baues, Dreckmann, Franjou
and Pirashvili-ring and the order we get.
As an application of our result, we count the
number of 'torsion free' indecomposable polynomial functors in
${\mathcal A}_{\hat\Z_p}^p$.

Our paper is organized as follows. In Section \ref{generalities}
we give the essential definitions and relate the different
concepts. In Section~\ref{knownfacts} we recall some of the most
important discoveries used in the sequel.
Section~\ref{recollement} describes the classical recollement diagrams
as well as the new one we have to use for ${\mathcal A}_{\F_p}$,
and we derive first consequences. The first main result is
proved in Section~\ref{modpcase}. We give the structure of
${\mathcal A}_{\F_p}^p$ there. Finally, in Section~\ref{lifttochar0}
we determine ${\mathcal A}_{\hat\Z_p}^p$ and prove the
second main theorem there.

\subsection*{Acknowledgement:} This research was
done during the years 2001 to 2003
as joint work with Steffen K\"onig. We presented the result on various occasions,
such as in Leicester in March 2003, Jena and Strasbourg in October 2003,
in Valenciennes in February 2004, in Bern in
February 2005, in Mainz in June 2005,
and in October 2005 in the ''s\'eminaire Chevalley'' Paris.
Recently we received numerous encouragements to
publish our manuscript.
Steffen K\"onig\footnote{email to the author from April 12, 2013}
wrote to me that he will not find the time to finish the paper,
and he gave me the autorisation to publish the paper alone. I wish to
thank Steffen K\"onig for having shared his insight with me, and
for allowing me to publish the paper.

I wish to thank the referee for many helpful comments and in particular for
indicating the treatment of Section~\ref{recollementforAnpsection}.

\section{Generalities on polynomial functors}
\label{generalities}

\subsection{Definitions}

Let $\mathcal A$ be a category with direct sums
and $\mathcal B$ be a category with direct sums and kernels. Then,
following Eilenberg and MacLane \cite{eilenberg}
define the cross effect $F^{(1)}$ of a functor $F:{\mathcal
A}\lra{\mathcal B}$ to be the bifunctor ${\mathcal A}\times
{\mathcal A}\lra{\mathcal B}$ defined on objects by
$F^{(1)}(V|W):=ker(F(V\oplus W)\lra F(V)\oplus F(W)),$ and on
morphism by the naturality of the construction. For $n\geq1$,
the $n$-th cross
effect $F^{(n)}$ is the cross effect of the $n-1$-st cross effect of
$F$, seen as functor in the first variable. Hence,
\begin{eqnarray*}
F^{(n)}(V|W|V_1|\dots|V_{n-1})&:=&ker\left(F^{(n-1)}(V\oplus
W|V_1|\dots|V_{n-1}) \lra\right.\\
&&
\left.\phantom{}F^{(n-1)}(V|V_1|\dots|V_{n-1})\oplus
F^{(n-1)}(W|V_1|\dots|V_{n-1})\right)
\end{eqnarray*}
for objects $V,W,V_1,\dots,V_{n-1}$ in $\mathcal A$.
Suppose in the sequel $F(0)=0$.
A functor $F$ is said to be polynomial of degree at most $n$ if
$F^{(n)}=0$ for an
$n\in\N$ (see Pirashvili \cite{pirashvilibourbaki}).
Given a commutative ring $R$, let $R-mod$
be the category of finitely generated $R$-modules, let $R-Mod$ be the category
of all $R$-modules, and let $R-free$ be the category
of finitely generated free $R$-modules.
Further, call ${\mathcal F}_R^n$
the category of polynomial
functors of degree at most $n$ from $R-mod$ to $R-mod$, and
${\mathcal F}_R:=\underrightarrow{\lim}_n\;{\mathcal F}_R^n$.
Moreover, let ${\mathcal A}_R^n$ be the category of polynomial
functors of degree at most $n$ from $\Z-free$ to $R-mod$, and
${\mathcal A}_R:=\underrightarrow{\lim}_n\;{\mathcal A}_R^n$. All these
categories are abelian.
Observe that additive functors are exactly the degree $1$
polynomial functors. The only degree $0$ polynomial functor is
the trivial functor due to our hypothesis that $F(0)=0$.

Friedlander and Suslin define in \cite[Definition 2.1]{friedlander}
the category of strict polynomial
functors ${\mathcal P}_k$ over a field $k$. A strict polynomial functor
$F$ is defined by associating to each finite dimensional $k$-vector space
a $k$-vector space $F(V)$ and to associate for any two finite
dimensional $k$-vector spaces $V$ and $W$ an element in
$$S^*\left(Hom_k(Hom_k(V,W),k)\right)\otimes Hom_k(F(V),F(W))$$
which in addition satisfy the usual compatibility
relations for compositions and the identity.
Each of these elements can be interpreted as mapping
$Hom_k(V,W)\lra Hom_k(F(V),F(W))$ by interpreting the formal
polynomial
$$S^*\left(Hom_k(Hom_k(V,W),k)\right)\otimes Hom_k(F(V),F(W))$$
as polynomial mapping, so that any strict polynomial functor induces a
polynomial functor $k-mod\lra k-Mod$.
Hence, we have a forgetful functor ${\mathcal P}_k\lra {\mathcal F}_k$.
It is shown in \cite{friedlander} that the
category of exact degree $n$ strict polynomial functors
${\mathcal P}_k^n$ from
$k-mod$ to $k-Mod$ for a field $k$ is equivalent to the category of modules
over the Schur algebra $S_k(n,n)$. Moreover,
${\mathcal P}_k=\bigoplus_n {\mathcal P}_k^n$.
Finally, the forgetful functor ${\mathcal P}_k\lra {\mathcal F}_k$
sends a strict polynomial functor of degree at most $d$ to a
polynomial functor of degree at most $d$
(cf \cite[Remark 4.1]{pirashvililecturenotes}).

A main theme in Section~\ref{generalities} in particular is the
question of base change. In Section~\ref{modpcase} we shall define
for any commutative rings $R$ and $S$ a particular degree $d$ polynomial functor
$$R[-]/(I^{d+1}):S-free\rightarrow R-Mod$$ which
assigns to every free $S$-module $V$ the quotient of the
group ring $R[V]$ of the additive group $V$ over $R$ by the $d+1$-th
power of the augmentation ideal $I$. Look at the case $S=\Z$.
A key idea, due to Pirashvili,
is to replace the morphisms $Hom_{\Z}(\Z^n,\Z^m)$ in the category of free abelian
groups by $R[Hom_{\Z}(\Z^n,\Z^m)]/(I^{d+1})$. This construction does
not allow base change. Strict polynomial functors basically replace $Hom_{\Z}(\Z^n,\Z^m)$
by $S^d(Hom_{\Z}(\Z^n,\Z^m))$, which does admit base change. Comparison between these two
constructions is the main theme in Section~\ref{generalities}.

\subsection{Functors with values in  characteristic $0$.}

\begin{Lemma}\label{isapolynomial} Let $R$ be an integral
domain of characteristic $0$.
If $F$ is a polynomial functor $F:\Z-free\lra R-mod$ of degree $d$,
then
$Hom_{\Z}(\Z^n,\Z^m)\stackrel{F}{\lra}Hom_{R}(F(\Z^n),F(\Z^m))$
is polynomial of degree $d$ in the $n\cdot m$ coordinate functions
$Hom_{\Z}(\Z^n,\Z^m)$.
\end{Lemma}

Proof. We shall show by induction on $n+d$ that for any $k$
homomorphisms $f_1,f_2,\dots,f_k\in Hom_{\Z}(\Z^n,\Z^m)$
and integers $\lambda_1,\lambda_2,\dots,\lambda_k\in\Z$ one gets
$F(\lambda_1f_1+\lambda_2f_2+\dots+\lambda_kf_k)$ is a degree $d$
polynomial in the variables $\lambda_1,\lambda_2,\dots,\lambda_k$.

If the degree of $F$ is $1$, there is nothing to show since then
the functor is linear.

Let $n>1$. Now, we know that
$$F(\Z^{n-1}\oplus\Z)=F(\Z^{n-1})\oplus F(\Z)\oplus
F^{(1)}(\Z^{n-1}|\Z)\;.$$
Hence, the restriction $f_i'$ of each of the $f_i$ to $\Z^{n-1}$
and the restriction $f_i''$ to the last component $\Z$ define
morphisms
$F(\sum_{i=1}^k\lambda_if_i'):F(\Z^{n-1})\lra F(\Z^m),$
$F(\sum_{i=1}^k\lambda_if_i''):F(\Z)\lra F(\Z^m)$
and
$F(\sum_{i=1}^k\lambda_if_i):F^{(1)}(\Z^{n-1}|\Z)\lra F(\Z^m)\;.$
In the first two cases, the dimension of the source space is
less than $n$, while the degree of $F$ is unchanged, whereas in
the third case the dimension of the source space is $n$, but the
degree of the functor is $d-1$. So, in any of these cases by the
induction hypothesis we can express $F(\sum_{i=1}^k\lambda_if_i)$
a polynomial of degree $n$ in the
variables $\lambda_1,\lambda_2,\dots,\lambda_k$.

We are left with the case $n=1$. The very same reduction applied to
the image and induction on $m+d$ implies that one can suppose that
$m=1$. But then, Eilenberg and MacLane
\cite[(8.3)]{eilenberg} show
that for $\lambda\in\Z$ one has
$$F(\lambda\cdot)=F((\lambda-1)\cdot)+F(1)+F^{(1)}(\lambda\cdot|1),$$
where $F(1)$ is the identity. Now,
$$F(\lambda\cdot)-F((\lambda-1)\cdot)=F(1)+F^{(1)}(\lambda\cdot|1),$$
where by the
induction hypothesis, since the degree of $F^{(1)}$ is less than
the degree of $F$, the right hand side $F(1)+F^{(1)}(\lambda\cdot|1)$ is
polynomial of degree $d-1$ in $\lambda$.

We shall now adapt an argument of Kuhn
\cite[Lemma 4.8]{Kuhn1}  to this slightly more general
situation. We claim that a function $f:\Z\lra R$ is a polynomial
if and only if some derivative $f^{(r)}$ vanishes, where
$f^{(r)}(n)=f^{(r-1)}(n)-f^{(r-1)}(n-1)$.

We assume for the moment that $R$ contains $\Q$.
Suppose $f$ is a polynomial. Then, it is clear that
$f^{(deg(f)+1)}=0$. Suppose to the contrary that $f^{(r)}=0$.
The polynomials
${X\choose k}:=\frac{X\cdot(X-1)\cdot\dots\cdot(X-k+1)}{k!}$
for $k\in\{0,1,\dots,d\}$
form an $R$-basis of the polynomials of degree at most $d$ in
$R[X]$, since $d!$ is invertible in $R$.
Moreover, ${X\choose k}-{X-1\choose k}={X\choose k-1}$.
Now, by induction, $f^{(s)}$ is a polynomial,
and hence a linear combination of polynomials ${X\choose k}$.
The relation ${X\choose k}-{X-1\choose k}={X\choose k-1}$ gives a
polynomial $h^{(s-1)}$ so that $(h^{(s-1)})^{(1)}=f^{(s)}.$
By induction on $n$, the values $g^{(1)}(n)=g(n)-g(n-1)$ determine
the values $g(n)$ up to the value of $g(0)$. Therefore, up to this
constant value, $f^{(s-1)}=h^{(s-1)}$.

Now, suppose $R$ an integral domain of characteristic $0$. Then,
since $R\subseteq frac(R)$, $f$ can be considered as being in
values $frac(R)$ which contains $\Q$.
This proves the claim.

Now, define $f(\lambda):=F(\lambda\cdot)$ and apply the claim
to conclude that $F(\lambda\cdot)$ is polynomial of
degree $d$.
\phantom{x}\hfill\dickebox

\bigskip

Let $K$ be a field.
A priori the category ${\mathcal A}_K$ is different from the
category ${\mathcal F}_K$. Nevertheless, in some cases we get one
inclusion.

\begin{Lemma}\label{inclusion} Let $K$ be either a prime field
of finite characteristic or let $K$ be a field of characteristic $0$.
Let $(\circ(K\otimes_{\Z}-))^*:{\mathcal F}_K\lra {\mathcal A}_K$
be the functor defined by
$(\circ(K\otimes_{\Z}-))^*(F):=F\circ(K\otimes_{\Z}-)$. Then,
$(\circ(K\otimes_{\Z}-))^*$ induces a fully faithful embedding
${\mathcal F}_K\hookrightarrow {\mathcal A}_K.$
\end{Lemma}

Proof:
The functor $\circ(-\otimes_\Z K):{\mathcal F}_K\lra {\mathcal A}_K$
induces for any two functors $F$ and $G$ in ${\mathcal F}_K$ a mapping
$$\varphi:Hom_{{\mathcal F}_K}(F,G)\lra
Hom_{{\mathcal A}_K}(F\circ(-\otimes_\Z K),G\circ(-\otimes_\Z K))\;.$$
We shall need to show that this mapping is an isomorphism.

\medskip

Injectivity:
Let $\eta_1$ and $\eta_2$ be two objects in $Hom_{{\mathcal F}_K}(F,G)$.
Suppose $\varphi(\eta_1)=\varphi(\eta_2)$. Observe that for any
$V\in\Z-free$ we have
$$(\varphi(\eta_1))(V)=
\eta_1(K\otimes_\Z V)\in Hom_K(F(K\otimes_\Z V),G(K\otimes_\Z V))$$
and likewise for $\eta_2$, satisfying that for any $V$ and $W$
and any $\rho\in Hom_\Z(V,W)$ one has
$(\eta_1(K\otimes_\Z V))\circ \hat G(\rho)=\hat F(\rho)\circ
(\eta_1(K\otimes_\Z W))$.

Since $\circ(-\otimes_\Z K)$ is dense,
$\eta_1$ and $\eta_2$ coincide on every object of $K-mod$, and so
$\eta_1=\eta_2$.

As a consequence, without any further hypothesis,
$$\varphi:Hom_{{\mathcal F}_K}(F,G)\hookrightarrow
Hom_{{\mathcal A}_K}(F\circ(-\otimes_\Z K),G\circ(-\otimes_\Z K))\;.$$

\medskip

Surjectivity:
Let
$\eta\in Hom_{{\mathcal A}_K}(F\circ(-\otimes_\Z K),G\circ(-\otimes_\Z K))$
be a natural transformation. We need to show that there is a natural
transformation $\eta'\in Hom_{{\mathcal F}_K}(F,G)$ so that
$\varphi(\eta')=\eta$.

In the case $K$ being a prime field of finite characteristic define
for any $V\in K-mod$ the mapping $\eta'(V):=\eta(P_V)$,
where $P_V$ is a fixed chosen projective cover of $V$, so that
$K\otimes_\Z P_V=V$.

In the case $K$ being of characteristic $0$, fix for any $K$-vector space
$V$ a free abelian subgroup $P_V$ so that $K\otimes_\Z P_V=V$.
Define $\eta'(V):=\eta(P_V)$.

We need to show that $\eta'$ is a natural transformation.

Let $\varphi\in Hom_K(V,W)$.

Consider first the case of $K$ being a prime field of
finite characteristic.
Since $P_V$ and $P_W$ are projective covers
of $V$ and $W$ as abelian groups, there is a
$\hat\varphi\in Hom_\Z(P_V,P_W)$ so that $K\otimes_\Z\hat\varphi=\varphi$
under the identification $K\otimes_\Z P_V=V$ and $K\otimes_\Z P_W=W$.
Since $\eta$ is a natural transformation,
$\hat G(\hat\varphi)\circ \eta(P_V)=\eta(P_W)\circ \hat F(\hat\varphi)$.
But, by definition, $\hat G(\hat\varphi)=G(\varphi)$ and
$\hat F(\hat\varphi)=F(\varphi)$,
as well as $\eta(P_W)=\eta'(W)$ and $\eta(P_V)=\eta'(V)$.
So, $\eta'$ is a natural transformation.

Suppose now that $K$ is a field of characteristic $0$.
Since ${\mathcal F}_K={\mathcal P}_K$ in this case, we know that
$F$ (and $G$ resp.) are polynomial laws transforming any linear
mapping $V\lra W$ into a linear mapping $F(V)\lra F(W)$ (and
$G(V)\lra G(W)$ resp.) which depends polynomially in the coefficients
of any matrix representation with respect to any fixed bases.
We know that for any $\Z$-linear mapping $\hat\varphi:P_V\lra P_W$
that the equation
$$(\ddag):G(K\otimes_\Z\varphi)\circ \eta(P_V)=
\hat G(\hat\varphi)\circ \eta(P_V)=
\eta(P_W)\circ \hat F(\hat\varphi)=\eta(P_W)\circ F(K\otimes_\Z\varphi).$$
Since this equation holds evaluated in infinitely coefficients, since
$\Z$ and $K$ are both infinite, this above equation $(\ddag)$ holds
{\em as polynomial} equation.

Therefore, the equation hold as well for $\varphi$, since there the
only difference is that the polynomials are evaluated not only on
integer coefficients, but also on coefficients in $K$. Since the
equation holds as polynomials, this equation holds true also
evaluated on $K$.

Therefore again $\eta'_W\circ F(\varphi)=G(\varphi)\circ \eta'_V$.
This proves that $\eta'$ is a natural transformation.
\phantom{x}\hfill\dickebox

\bigskip

We are now concerned with the question when a polynomial functor
$\Z-free\lra R-mod$  can be extended to a polynomial functor
$R-mod\lra R-mod$ by composing with the 'extending scalars' functor
$\Z-free\stackrel{R\otimes_{\Z}-}{\lra}R-mod$. In other words
we study the question when $(\circ(K\otimes_{\Z}-))^*$ is an
equivalence ${\mathcal A}_R\simeq {\mathcal F}_R$. In order to
prove this, by Lemma~\ref{inclusion}, one needs to show that
$(\circ(K\otimes_{\Z}-))^*$ is dense as well.

We have to deal with mainly two cases: the case of $R$ being a
field of characteristic $0$ and the case of $R$ being a field of
characteristic $p$. We shall see that fields with characteristic
$0$ behave more like characteristic $\infty$. The remarks at the beginning
of this section on base change properties
are particularly visible in the proof of
the following lemma.

\begin{Lemma}\label{char0} Let $R$ be a field
of characteristic $0$ and let $F:\Z-free\lra R-mod$ be a polynomial
functor of degree $d$. Then, $F$ extends to a polynomial functor
$\hat{F}:R-free\lra R-mod$ so that $\hat{F}\circ
(R\otimes_{\Z}-)=F$. In particular, $(\circ
(R\otimes_{\Z}-))^*$ induces an equivalence
${\mathcal A}_R\simeq {\mathcal F}_R$.
\end{Lemma}

Proof.
By Lemma~\ref{inclusion} we know that ${\mathcal F}_K^d\hookrightarrow
{\mathcal A}_K^d$.

We have to show that this embedding is dense.
Let $F$ be a degree $d$
polynomial functor in ${\mathcal A}_K^d$. By Lemma~\ref{isapolynomial}
we know that for any $n$ and $m$ the functor
$F$ induces a degree $d$ polynomial mapping with coefficients in $K$
in the coordinate functions of matrices in $Hom_\Z(\Z^n,\Z^m)$. Moreover,
$F$ is a functor, that is $F(\alpha \beta)=F(\alpha)F(\beta)$
and $F(id_A)=id_A$ for any free abelian group $A$ and any two composable
morphisms of abelian groups $\alpha$ and $\beta$.

Let $\alpha:\Z^n\lra \Z^m$ and $\beta:\Z^k\lra \Z^n$. The equation
$F(\alpha \beta)=F(\alpha)F(\beta)$ translates into an equation
between the evaluation of the corresponding polynomials in each degree.
Since $\Z$ and $K$ are of characteristic $0$,
the polynomial equation holds
if evaluated on infinitely many values, and so the polynomial equations
actually holds {\em as polynomials}.
Friedlander and Suslin remark in
\cite[remark after Definition 2.1]{friedlander} that this is
actually equivalent to saying that $F$ actually is a strict polynomial
functor $\hat G\in {\mathcal P}_K$  of degree $d$.
Now, for $K$ infinite, ${\mathcal P}_K^{\leq d}\simeq {\mathcal F}_K^d$.
(\cite{friedlander}). So, actually $\hat G$, and in turn $F$
is a degree $d$ polynomial functor in ${\mathcal F}_K^d$.
This proves the lemma.\hfill\dickebox

\subsection{Functors with values in fields of finite
characteristic.}

\begin{Lemma}\label{welldefined}
Let $\F$ be a field of characteristic $p$ and
let $F:\Z-free\lra \F-mod$ be polynomial
functor  of degree at most $p-1$ which preserves the initial
object, i.e. $F(0)=0$. Then, for any homomorphism $\alpha:\Z^n\lra
\Z^m$ one gets $F(p\cdot \alpha)=0$.
\end{Lemma}

Proof. Let $M=\Z^n$ and $N=\Z^m$.
We write $p\cdot \alpha$ in the diagram
$$\begin{array}{ccccc}
M&\stackrel{\delta_M}{\lra}&M\oplus M\oplus\dots\oplus
M&\stackrel{\sigma_M}{\lra}&M\\
\mbox{\scriptsize $\alpha$}\dar\phantom{\mbox{\scriptsize $\alpha$}}&&\phantom{\mbox{\scriptsize $\alpha\oplus\alpha\oplus\dots\oplus\alpha$}}\dar\mbox{\scriptsize $\alpha\oplus\alpha\oplus\dots\oplus\alpha$}&&\phantom{\mbox{\scriptsize $\alpha$}}\dar\mbox{\scriptsize $\alpha$}\\
N&\stackrel{\delta_N}{\lra}&N\oplus N\oplus\dots\oplus N&
\stackrel{\sigma_N}{\lra}&N
\end{array}
$$
where $\delta$ is the diagonal mapping and $\sigma$ is the
summation mapping, on $M$ or on $N$ respectively.

Denote by $F^{(i)}$ the $i^{th}$ cross effect of the functor $F$.
Since $F$ is polynomial  of degree at most $p$,
$$
F(A_1\oplus A_2\oplus\dots\oplus A_p)=
\bigoplus_{i=1}^p\bigoplus_{j_1<\dots<j_i}F^{(i-1)}(A_{j_1}|\dots|A_{j_i})
$$
for $p$ abelian groups $A_1,\dots,A_p$, and this decomposition
is functorial with respect to these groups. Moreover, since
$deg(F)=p$, and since $F(0)=0$, one gets $F^{(p-1)}=0$.

Since $F(p\cdot\alpha)=F(\alpha\circ\sigma_M\circ\delta_M)=
F(\alpha)\circ F(\sigma_M\circ\delta_M)$, it is necessary and
sufficient to show that $F(\sigma_M\circ\delta_M)=0$.
Now,
$$
F(M\oplus M\oplus\dots\oplus M)=
\bigoplus_{i=1}^p\bigoplus_{j_1<\dots<j_i}F^{(i-1)}(M_{j_1}|\dots|M_{j_i})
$$
where $M_{l}=M$ is the copy of $M$ in the $l$-th position of
$\oplus_{i=1}^pM$, for all $l\leq p$.
The mapping $F(\sigma_M\circ\delta_M)=F(\sigma_M)\circ
F(\delta_M)$ factors as a sum
$\sum_{i=1}^pF(\sigma_M)\circ \iota_i\circ \pi_i\circ F(\delta_M)$
where $\iota_i$ is the embedding of
$\bigoplus_{j_1<\dots<j_i}F^{(i-1)}(M_{j_1}|\dots|M_{j_i})$ into
$F(M\oplus M\oplus\dots\oplus M)$ and $\pi_i$ is the projection of
$F(M\oplus M\oplus\dots\oplus M)$ to this direct factor.
But now, for $i<p$,
$F(\sigma_M)|_{(\bigoplus_{j_1<\dots<j_i}F^{(i-1)}(M_{j_1}|\dots|M_{j_i}))}
\circ\pi_i\circ F(\delta_M)$ is a sum of $p$ identical mappings,
which sum up to $0$ in characteristic $p$.
Hence,
$F(\sigma_M\circ\delta_M)=
F(\sigma_M)\circ F^{(p-1)}(\alpha|\alpha|\dots|\alpha)\circ F(\delta_M)=0$
using that $F^{(p-1)}=0$.
\phantom{x}\hfill\dickebox

\begin{Cor}\label{modpwelldefined}
Let $\F$ be a field of characteristic $p$ and
let $F:\Z-free\lra \F-mod$ be polynomial
functor  of degree at most $p-1$ which preserves the initial
object, i.e. $F(0)=0$. Then, for any two homomorphisms $\alpha:\Z^n\lra
\Z^m$ and $\beta:\Z^n\lra\Z^m$ so that
$\alpha-\beta\in p\cdot Hom_{\Z}(\Z^n,\Z^m)$,
one gets $F(\alpha)=F(\beta)$.
\end{Cor}

Proof. This is a consequence of the previous lemma and
\cite[p. 76, formula 8.5]{eilenberg} and \cite[Theorem 9.3]{eilenberg}.
Indeed,
$F(\sum_{n=1}^{p+1}\rho_n)=
\sum_{n=1}^{p+1}\sum_{i_1<\dots<i_n}F^{(n-1)}(\rho_{i_1}|\dots|\rho_{i_n})$
implies
$$
F(\alpha+p\cdot\gamma)=F(\alpha)+F^{(p-1)}(\underbrace{\gamma|\gamma|\dots
|\gamma}_{p\mbox{ \scriptsize factors}})+
F^{(p)}(\alpha|\underbrace{\gamma|\gamma|\dots
|\gamma}_{p\mbox{ \scriptsize factors}})
=F(\alpha)+F(p\gamma)=F(\alpha)
$$
by Lemma~\ref{welldefined}.
\hfill\dickebox

\bigskip

\begin{Lemma}\label{charp}
Let $\F$ be the prime field of characteristic $p$ and
let $F:\Z-free\lra \F-mod$ be a polynomial functor of degree
less or equal to $p-1$.
Then, $F$ factors through the functor
$\F\otimes_{\Z}-:\Z-free\lra\F-mod$. Moreover, if
$F=F'\circ(\F\otimes_{\Z}-)$, then $F$ is polynomial of degree $m$
if and only if $F'$ is polynomial of degree $m$. Hence,
$(\circ (\F\otimes_{\Z}-))^*$ induces an equivalence
${\mathcal A}_{\F}^{\leq p-1}\simeq {\mathcal F}_{\F}^{\leq p-1}$.
\end{Lemma}

Proof. Since $(\F\otimes_{\Z}-)^*:{\mathcal F}_{\F}\hookrightarrow
{\mathcal A}_{\F}$ is a fully faithful embedding, we need to show
that $(\F\otimes_{\Z}-)^*$ is dense. So, given a functor
$F:\Z-free\lra  \F-mod$. One has to show that there is a functor
$\hat F:\F-mod\lra \F-mod$ with $\hat F\circ (\F\otimes_{\Z}-)=F$.

For any $V\in \F-mod$ choose $P_V$ a projective cover as abelian group.
Then, $P_V$ is in $\Z$-free. By the universal property of projective
covers one has for any $\alpha\in Hom_\F(V,W)$ a (non-unique)
$\hat\alpha\in Hom_\Z(P_V,P_W)$ so that
$$\begin{array}{ccc}P_V&\stackrel{\hat\alpha}{\lra}&P_W\\
\dar&&\dar\\ V&\stackrel{\alpha}{\lra}&W\end{array}$$ is commutative.
Put
$$\hat F(V):=F(P_V)\;\mbox{ and }\hat F(\alpha):=F(\hat\alpha)\;.$$
We need to show that this gives a functor $\hat F:\F-mod\lra\F-mod$.
Let $\hat\alpha$ and $\hat\alpha'$
be two different lifts of $\alpha:V\lra W$,
then $\hat\alpha-\hat\alpha'$ lifts the $0$-mapping, and so
$\hat\alpha-\hat\alpha'\in p\cdot Hom_\Z(P_V,P_W).$
Corollary~\ref{modpwelldefined} implies that
$\hat F(\hat\alpha)=\hat F(\hat\alpha')$.
Using Corollary~\ref{modpwelldefined} again one gets that
$\hat F(id_V)=id_{\hat F(V)}$ since $id_{P_V}$ is a lift of $id_V$.
Moreover, let $\alpha:U\lra V$ and $\beta:V\lra W$, then
choosing lifts $\hat\alpha:P_U\lra P_V$ and $\beta:P_V\lra P_W$, one gets
$\widehat{\alpha\beta}-\hat\alpha\hat\beta$ lifts the $0$-mapping.
So, $\widehat{\alpha\beta}-\hat\alpha\hat\beta\in p\cdot Hom_\Z(P_U,P_W)$
and again by Corollary \ref{modpwelldefined}
one has $F(\widehat{\alpha\beta})=F(\hat\alpha\hat\beta)$.
\phantom{x}\hfill\dickebox

\section{A review on polynomial functors and functor cohomology}

\label{knownfacts}

\subsection{Polynomial functors are modules}
\label{bdfp}

Let $R$ be a commutative ring.
We know that by a result of Baues, Dreckmann, Franjou and Pirashvili
\cite{BDFP} a polynomial functor of degree at most $n$
from free abelian groups to $R$-modules is defined by giving
$R$-modules $F_m(\Z|\Z|\dots|\Z)$ for all $m\leq n$ and  mappings
$$h_k^m: F_m(\Z|\Z|\dots|\Z)\lra F_{m+1}(\Z|\Z|\dots|\Z)$$
for $k\leq m\leq n-1$ and
$$p_k^{m+1}: F_{m+1}(\Z|\Z|\dots|\Z)\lra F_{m}(\Z|\Z|\dots|\Z)$$
for $k\leq m\leq n-1$ satisfying the  relations~:

$$(*)\;\;h_j^mp_i^m=
\left\{\begin{array}{ll}
p_i^{m+1}h_{j+1}^{m+1}&\mbox{ for }j<i\\
p_{i+1}^{m+1}h_{j}^{m+1}&\mbox{ for }j>i\\
1+t_i^m+p_i^{m+1}h_{i+1}^{m+1}+
p_{i+1}^{m+1}h_i^{m+1}+p_{i+1}^{m+1}t_i^{m+1}h_{i+1}^{m+1}\\
+p_{i}^{m+1}t_{i+1}^{m+1}h_{i}^{m+1}+
p_{i+1}^{m+1}p_{i}^{m+2}t_{i+1}^{m+2}h_{i}^{m+2}h_{i+1}^{m+1}
&\mbox{ for }j=i\end{array}\right\}
$$

Define the algebra $\Gamma^n_R$ over $R$
by a quiver with $n$ vertices
$F_m(\Z|\Z|\dots|\Z)$ for any $m$ with $1\leq m\leq n$ and arrows
$h_k^m: F_m(\Z|\Z|\dots|\Z)\lra F_{m+1}(\Z|\Z|\dots|\Z)$ and
$p_k^{m+1}: F_{m+1}(\Z|\Z|\dots|\Z)\lra F_{m}(\Z|\Z|\dots|\Z)$
subject to the relations $(*)$.
Observe that the relations above do not form a set of admissible
relations. The third relation though should be read as the defining
equation for the symbols $t_i^n$, and this way the relations
$(*)$ is a set of admissible relations.
The result \cite{BDFP} of Baues et alii implies that
$\Gamma^n_R-mod$ is equivalent to the category of polynomial functors
of degree at most $n$.

Since the ring homomorphism $\hat Z_p\lra\F_p$ induces an
embedding $\F_p-mod\lra \hat\Z_p-mod$, we get an induced embedding
$\Gamma^n_{\F_p}-mod\lra \Gamma^n_{\hat\Z_p}-mod$ which is also induced
by the surjective ring homomorphism
$\Gamma^n_{\hat\Z_p}\lra \Gamma^n_{\F_p}$.

\begin{Rem}
We should mention that the description of \cite{BDFP}
was recently generalised by Hartl, Pirashvili and Vespa \cite{HPV}
to show an isomorphism of the category
of functors $\Z-free\lra \Z-Mod$
with the category of what they call pseudo-Mackey functors.
\end{Rem}

\subsection{Some facts on functor cohomology}
\label{functorcohomology}

We shall give some facts that we will need from Franjou, Friedlander,
Scorichenko and Suslin \cite{FFSS}. Basically, these results reduce
the computation of extension groups between polynomial functors to
questions between extension groups between strict polynomial functors.

Let $k$ be a commutative ring.
Let ${\cal F}(k)$ be the category of functors $k-mod\lra k-mod$
and let ${\cal P}(k)$ the category of strict polynomial functors
between $k$-modules.
Let ${\cal F}^n={\cal F}^n(k)$ be the category of degree $n$
polynomial functors from $k-mod$ to $k-mod$.
If $k=\F_q$ for $q=p^s$, the field with $q$ elements, we write
${\cal F}^n(k)={\cal F}^n(q)$.
In this case, for any strict polynomial functor $P$ in
${\cal P}(k)$ let $P^{(m)}$
be the functor twisted by the Frobenius endomorphism defined by
$\F_q\ni x\mapsto x^{p^m}\in\F_q$.

\begin{Theorem}\cite{FFSS}\label{FFSSresult}
Given any two homogeneous strict polynomial functors
$P$ and $Q$ between $\F_q$-vector spaces for $q=p^s$.
If the degrees of $P$ and $Q$ are different and strictly smaller than
$q$, then $Ext^*_{{\cal F}(k)}(P,Q)=0$.
Moreover, if the degree of $P$ and $Q$ coincide, then
$${\lim_{\longrightarrow}}_mExt^*_{{\cal P}(k)}(P^{(m)},Q^{(m)})\simeq
Ext^*_{{\cal F}(k)}(P,Q)$$
\end{Theorem}

Frobenius twisting decreases the 'degree of homological
triviality' as is shown in a result of H.H. Andersen.

\begin{Prop}(H. H. Andersen;
see \label{Andersen}\cite[Corollary 1.3]{FFSS})
For two homogeneous strict polynomial functors
$P$ and $Q$ between $\F_q$-vector spaces of the same degree,
for $q=p^s$ and for $m\in\N\cup\{0\}$ we get
$$Ext^*_{{\cal P}(k)}(P^{(m)},Q^{(m)})\leq
Ext^*_{{\cal P}(k)}(P^{(m+1)},Q^{(m+1)})$$
\end{Prop}

The first of the two statements in Theorem \ref{FFSSresult}
actually is due to Kuhn:

\begin{Lemma}\cite{Kuhn1} \label{Kuhnlemma}
Any functor $F\in{\cal F}(q)$ decomposes into a
direct sum $F=\oplus_{i=0}^{q-1}F_i$ where
$F_i(V):=\{x\in F(V)|\;F(\lambda\cdot)(x)=\lambda^i\cdot x\;\;
\forall\lambda\in\F_q\}$. This induces a decomposition
of the category of functors between $\F_q$-vector spaces
${\mathcal F}(q)=\prod_{i=0}^{q-1}{\cal F}(q)_i$.
\end{Lemma}

Finally, a result due to Kuhn will be essential in the sequel.

\begin{Theorem}(N. Kuhn \cite{KuhnII,kuhnsurvey})
\label{uniserials}\label{Kuhnresult}
The injective envelope $I_{\F_p}$ of the trivial module in the category of
analytic functors ${\cal F}^{\omega}(\F_p)$ from finite
dimensional $\F_p$-vector spaces to $\F_p$-vector spaces is
uniserial and the only composition factors in ${\cal F}^p(\F_p)$
are the two composition factors of $soc_2(I_{\F_p}))$, where, as usual,
$soc_2$ denotes the second layer in the socle series \cite[Definition 1.2.1]{Benson}.
\end{Theorem}

The next result of Franjou, Lannes and Schwartz implies that the
categories ${\cal A}^{p}_{\F_p}$ and ${\cal F}^{p}_{\F_p}$ are
different.

\begin{Theorem}(Franjou-Lannes-Schwartz \cite{FLS}; Franjou-Pirashvili
\cite{FP})\label{FLS}
$$Ext_{{\cal A}^{p}_{\F_p}}(id,\F_p\otimes id)\simeq
\F_p[e_1,e_2,\dots]/(e_h^p; h\geq 0)\otimes \Lambda(\xi_1)$$
where
$\Lambda(\xi_1)$ is the exterior algebra in one variable with
generator in degree $2p-1$ and $e_h$ are generators in degree
$2p^h$. Moreover,
$$Ext_{{\cal F}^{p}_{\F_p}}(id, id)\simeq
\F_p[e_0,e_1,\dots]/(e_h^p; h\geq 0).$$
\end{Theorem}

Remark that $\F_p\otimes_\Z id=id$ as functors on the category $\F_p-mod$.
As a consequence, $Ext^2_{{\mathcal F}^p_{\F_p}}(id,id)\neq 0$.
Indeed, by \cite[7.3]{FLS} the following four term sequence is a
non zero element:
$$0\lra id\lra S_p\lra S^p\lra id\lra 0;$$
where $S_p$ is the degree $p$ homogeneous part of the
coinvariants under the $\Sn_p$ action on the tensor algebra,
and where $S^p$ is the degree $p$ homogeneous part of the
invariants of the tensor algebra.

\section{On recollement diagrams}
\label{recollement}

We remind the reader to the notion of a recollement diagram.
A recollement diagram is given by
three
categories ${\mathcal A},{\mathcal B}$ and ${\mathcal C}$ with
functors
$$
{\mathcal A}
\begin{array}c\stackrel{q}{\leftarrow}\\
\stackrel{i}\rightarrow\\ \stackrel{p}\leftarrow \end{array}
{\mathcal B}
\begin{array}c\stackrel{l}{\leftarrow}\\
\stackrel{e}\rightarrow\\ \stackrel{r}\leftarrow \end{array}
{\mathcal C}
$$
so that
\begin{enumerate}
\item $(l,e)$ and $(e,r)$ are adjoint pairs.
\item $(q,i)$ and $(i,p)$ are adjoint pairs.
\item $i$ is a full embedding and $e(B')=0\Leftrightarrow B'\simeq
i(A')$ for an $A'\in\mathcal A$.
\item the adjointness morphisms $e\circ r\lra id_{\mathcal C}$ and
$id_{\mathcal C}\lra e\circ l$ are isomorphisms.
\end{enumerate}
We denote a recollement diagram as above by
$\left({\mathcal A},{\mathcal B},{\mathcal C},(e,l,r),(i,q,p)\right)$.

\medskip

We shall give a result which is a special case of a recent result of
Chrysostomos Psaroudakis and Jorge Vitoria \cite{Chrysostomos}.
We shall give our original proof below, since in our special case the proof
is much easier than the proof for the
general statement from \cite{Chrysostomos}.

\begin{Prop}\label{recollementofmodules}
Let $\left({\mathcal A},{\mathcal B},{\mathcal C},(e,l,r),(i,q,p)\right)$
be a recollement diagram.
Suppose ${\mathcal B}\simeq B-mod$ and ${\mathcal C}=C-mod$
are module categories and that ${\mathcal B}$
satisfies the Krull-Schmidt theorem on
projective modules. Suppose that $e$ is representable. Then,
${\mathcal A}=A-mod$ again is a module category, there is
an idempotent $e'$ in $B$ so that $C$ is Morita equivalent to
$e'Be'$,  and $A$ is Morita equivalent to $B/Be'B$.
\end{Prop}

Proof.
Since $e$ has a left and a right adjoint, $e$ is exact. Therefore,
$e=Hom_{\mathcal B}(P,-)$ where $P$ is a projective
object. Since $l$ and $r$ are left and right adjoints to $e$, we
get that $l=P\otimes_{End_B(P)}-$, that $r=Hom_{End_B(P)}(Hom_B(P,B),-)$
and that $C\simeq End_B(P)$. Since ${\mathcal B}$ is a Krull-Schmidt
category, then up to Morita equivalence, one can choose $P=Be'$ for
an idempotent $e'^2=e'\in B$ and we get $C$ is Morita equivalent to
$e'Be'$.

The third condition in the definition of a recollement diagram
implies that ${\mathcal A}$ can be identified with those $B$-modules
$M$ for which $e'M=0$. Hence,
$${\mathcal A}\simeq\{M\in B-mod|\;e'M=0\}\simeq B/Be'B-mod.$$
This proves the proposition.
\phantom{x}\hfill\dickebox

\bigskip

Another important observation is
that, by the adjointness properties,
$l$ maps projective object in ${\mathcal C}$ to
projective objects in ${\mathcal B}$, and that
$r$ maps injective objects in ${\mathcal C}$ to
injective objects in ${\mathcal B}$.

\subsection{Analyzing Schwartz' recollement for polynomial functors}
\label{recollerpolys}

Let $q=p^s$ for a prime $p$. Then, using the notation of
Section~\ref{functorcohomology},
following Kuhn \cite[Theorem 1.3]{stratifying} or
Schwartz \cite[\S 5.5]{unstable} we have a recollement diagram
$$
{\cal F}^{n-1}(q){
\begin{array}c\leftarrow\\ \rightarrow\\ \leftarrow \end{array}}
{\cal F}^{n}(q){
\begin{array}c\leftarrow\\ \rightarrow\\ \leftarrow \end{array}}
\prod_{n(\lambda)=n}{\mathbb F}_q\Sn_{\lambda}-mod
$$
where $\Sn_{\lambda}=\Sn_{\lambda_1}\times\dots\times\Sn_{\lambda_{s-1}}$
where $\Sn_k$ is the symmetric group on $k$ elements and where
$n(\lambda):=\lambda_0+\dots+\lambda_{s-1}$.
Moreover, the functor ${\cal F}^{n}(q)\lra \prod_{n(\lambda)=n}{\mathbb
F}_q{\Sn_{\lambda}}-mod$ is representable by $id^{\lambda}$
and for a partition $\lambda=(\lambda_0\geq\dots\geq\lambda_{s-1})$, we set
$id^{\lambda}:=\oplus_{j=1}^{s-1} id^{\otimes\lambda_j}$.

\begin{Rem}
\label{nkleinerp}
Hence, in case $s=1$ and $n<p$, the recollement becomes
$$
{\cal F}^{n-1}(p){
\begin{array}c\leftarrow\\ \rightarrow\\ \leftarrow \end{array}}
{\cal F}^{n}(p){
\begin{array}c\leftarrow\\ \rightarrow\\ \leftarrow \end{array}}
\prod_{ \mbox{\scriptsize  partitions of }n}
({\mathbb F}_p-mod)
$$
since $\F_p\Sn_n$ is semisimple, since $\F_p$ is a splitting field,
and therefore its module category is equivalent to a direct product of copies
${\mathbb F}_p-mod$.
\end{Rem}

In case $s=1$ and $n=p$, the recollement becomes
$$
{\cal F}^{p-1}(p){
\begin{array}c\leftarrow\\ \rightarrow\\ \leftarrow \end{array}}
{\cal F}^{p}(p){
\begin{array}c\leftarrow\\ \rightarrow\\ \leftarrow \end{array}}
{\mathbb F}_p{\Sn_p}-mod
$$

\bigskip

We have an immediate consequence of the above result.

\begin{Lemma}\label{simples}
For any simple
polynomial functor $F$ there is a strict polynomial functor
$\hat{F}$ such that the forgetful functor, which assigns to every
strict polynomial functor its polynomial functor by evaluating the
polynomial as mapping, maps $\hat F$ to $F$.
\end{Lemma}

Proof. This is done by induction on the degree.
The simple objects in ${\cal F}^d(q)$ are in bijection with the union of
the simple objects in ${\cal F}^{d-1}(q)$ and the simple objects in
$\F_q\Sn_d-mod$. Now, any simple $\F_q\Sn_d$-module
is image of a simple module of the Schur algebra $S_{\F_q}(d,d)$
under the Schur functor as is a
classical fact. Since the category of degree $d$ strict polynomial
functors is equivalent to the category of modules over the Schur algebra
$S(d,d)$,  the simple objects in ${\cal F}^d(q)$ of degree $d$ are images of
a strict polynomial functor. By induction, the simple objects of degree
less than $d$ are images of strict polynomial functors.
\phantom{x}\hfill\dickebox

\subsection{Recollement for ${\mathcal A}_{\F_p}$}

\label{recollementforAnpsection}

For the category ${\mathcal A}_{\F_p}$ we get a similar
recollement diagram.

\begin{Rem}
The treatment we give for Section~\ref{recollementforAnpsection}
was suggested by the referee, and follows paths which are not
easily documented. They appeared in the special case of degree $2$ in introductory
remarks in \cite{Franjouext,FLS}. An alternative proof can be obtained, and actually this
was our initial approach,
using Piriou's thesis~\cite{piriouthese,PS} along its first chapter.
I am very grateful to the referee for this useful hint.

The recollement diagram was proved in a more general situation in the meantime by
Djament and Vespa~\cite[Theorem 2.2]{DjamentVespa2}, so that the
result Proposition~\ref{recollementA}
is a special case of \cite[Theorem 2.2]{DjamentVespa2}.
\end{Rem}

We shall mainly work with properties of adjoint functors.
Lemma~\ref{adjointsonfunctorlevel} below seems to be a well-known result.
It appears maybe for the first time in the proof of
\cite[Lemma 0.4]{FLS} without further reference. In order to keep the presentation as self-contained as possible, we shall provide a short proof.

\begin{Lemma} \label{adjointsonfunctorlevel}
Let $\mathcal C$, ${\mathcal C}'$ and $\mathcal D$ be categories, let
$L:{\mathcal C}\lra {\mathcal C}'$ be a functor admitting a right adjoint $R$, and let
$\hat L$ be the induced functor
$Funct({\mathcal C}',{\mathcal D})\lra Funct({\mathcal C},{\mathcal D})$ on the functor categories
given by precomposition with $L$, and let likewise
$\hat R$ be the induced functor
$Funct({\mathcal C},{\mathcal D})\lra Funct({\mathcal C}',{\mathcal D})$ given by precomposition with $R$.
Then $\hat R$ is left adjoint to $\hat L$.
\end{Lemma}

Proof.
We get natural transformations
$$id_{\mathcal C}\stackrel{\eta}{\lra} RL\mbox{ and }
LR\stackrel{\epsilon}{\lra} id_{{\mathcal C}'}$$
so that the compositions $$L\lra LRL\lra L \mbox{ and }R\lra RLR\lra R$$
are the identity on the respective ctegories (cf e.g. Maclane \cite[IV.1.Theorem 2]{Maclane}).

Now, given any functor $X:{\mathcal C}\lra {\mathcal D}$, then
we obtain a natural transformation
$$id_{Funct({\mathcal C},{\mathcal D})}\lra (\hat L\hat R)$$
induced from $\eta$
by $$X\stackrel{X(\eta)}{\lra} XRL=(\hat L\hat R)(X)$$ and
likewise a natural transformation
$$(\hat R\hat L)\lra id_{Funct({\mathcal C}',{\mathcal D})}.$$

Since the compositions $L\lra LRL\lra L$ and $R\lra RLR\lra R$
are the identity, this holds as well for the compositions
$\hat L\lra \hat L\hat R\hat L\lra \hat L$ and
$\hat R\lra \hat R\hat L\hat R\lra \hat R$.
Again by \cite[IV.1.Theorem 2]{Maclane} we obtain the statement.\hfill \dickebox

\bigskip

We shall need a very simple observation.

\begin{Lemma}\label{degreeunderext}
Let $A$ and $B$ be two polynomial functors of degree at most $n$
and let $$0\lra A\lra C\lra B\lra 0$$ be an exact sequence of functors.
Then, the degree of $C$ is at most $n$ as well.
Moreover, taking cross effects is exact.
\end{Lemma}

Proof: We get a commutative diagram
$$
\begin{array}{ccccccccc}
&&0&&0&&0\\
&&\dar&&\dar&&\dar\\
0&\lra&A^{(2)}(U|V)&\lra&C^{(2)}(U|V)&\lra&B^{(2)}(U|V)\\
&&\dar&&\dar&&\dar\\
0&\lra&A(U\oplus V)&\lra&C(U\oplus V)&\lra&B(U\oplus V)&\lra&0\\
&&\dar&&\dar&&\dar\\
0&\lra&A(U)\oplus A(V)&\lra&C(U)\oplus C(V)&\lra&
B(U)\oplus B(V)&\lra&0\\
&&\dar&&\dar&&\dar\\
&&0&&0&&0
\end{array}
$$
and the snake lemma implies that
$$0\lra A^{(2)}(U|V)\lra C^{(2)}(U|V) \lra B^{(2)}(U|V)\lra 0$$
is exact. Induction on the degree gives the result.
\hfill\dickebox

\bigskip

We shall now prove a general statement on functor categories which is an adaption from
Franjou~\cite[Section 1]{Franjouext}.

\begin{Lemma} \label{macdonaldlinearisation}
Let $R$ be a commutative ring and let
$G$ is an object in ${\mathcal A}^n_{R}$. Then we get
$Hom_{{\cal A}_{R}}((R\otimes id^{\otimes n},G)=G^{(n-1)}$.
\end{Lemma}

Proof.
Let
$$n-\Z-free=\underbrace{(\Z-free)\times\dots\times(\Z-free)}_{n\mbox{ \scriptsize $copies$}}$$
be the category with objects
$F_1\times\dots\times F_n$ for $F_i$ being an object in $\Z-free$ for each $i\in\{1,\dots,n\}$ and morphisms being
$$Hom_{n-\Z-free}(F_1\times\dots\times F_n,G_1\times\dots\times G_n):=
Hom_\Z(F_1,G_1)\times \dots\times Hom_\Z(F_n,G_n),$$
for all objects $F_i,G_j$ of $\Z-free$.
Composition of morphisms is given by composition of mappings in $\Z-free$.
We shall define $$\Pi_n:n-\Z-free\lra \Z-free$$ and $$\Delta_n:\Z-free\lra n-\Z-free$$
by $\Pi_n(F_1,\dots,F_n):=F_1\oplus\dots\oplus F_n$ and
$\Delta_n(F):=(F,\dots,F)$. Then $\Pi_n$ is left and right adjoint to $\Delta_n$
as is easily seen, almost by definition.

Now, we consider the category $n-{\cal A}_{R}$ of functors $n-\Z-free\lra R-Mod$
with morphisms being natural transformations.
We get functors
\begin{eqnarray*}
n-{\cal A}_{R}&\stackrel{\Delta^n}{\lra}&{\cal A}_{R}\\
F&\mapsto&F\circ\Delta_n
\end{eqnarray*}
and
\begin{eqnarray*}
{\cal A}_{R}&\stackrel{\Pi^n}{\lra}&n-{\cal A}_{R}\\
F&\mapsto&F\circ\Pi_n
\end{eqnarray*}

\medskip

By Lemma~\ref{adjointsonfunctorlevel}, we  see that
the functor $\Pi^n$ is left and right adjoint to $\Delta^n$.
Put
\begin{eqnarray*}
n-\Z-free&\stackrel{\boxtimes^n}{\lra}&\Z-free\\
(M_1,\dots,M_n)&\mapsto&M_1\otimes\dots\otimes M_n
\end{eqnarray*}
and obtain $id^{\otimes n}=\boxtimes^n\circ\Delta_n$ and
$R\otimes id^{\otimes n}=(R\otimes\boxtimes^n)\circ\Delta_n$.
But then
\begin{eqnarray*}
Hom_{{\cal A}_{R}}((R\otimes id^{\otimes n},G)&=&
Hom_{{\cal A}_{R}}((R\otimes\boxtimes^n)\circ\Delta_n,G)\\
&=&
Hom_{n-{\cal A}_{R}}((R\otimes\boxtimes^n),G\circ\Pi_n)
\end{eqnarray*}

\medskip

Let $I:=\{i_1,\dots,i_m\}\subseteq \{1,\dots,n\}$ be an $m$-element subset
of $\{1,\dots,n\}$, and suppose $i_1<\dots<i_m$. Then we may consider
$\pi_I^n:n-\Z-free\lra m-\Z-free$ to be the functor given by
$$\pi_I^n(M_1,\dots,M_n):=(M_{i_1},\dots,M_{i_m})$$
and the functor
$\iota_I^n:m-\Z-free\lra n-\Z-free$ given by injection into the corresponding
coordinates so that $\pi_I^n\circ\iota_I^n=id$. Then, again by definition,
$\iota_I^n$ is left and right adjoint to $\pi_I^n$.
Using Lemma~\ref{adjointsonfunctorlevel} we get that the functors
$$\widehat{\iota_I^n}:n-{\mathcal A}_{R}\lra m-{\mathcal A}_{R}$$
and
$$\widehat{\pi_I^n}:m-{\mathcal A}_{R}\lra n-{\mathcal A}_{R}$$
obtained by pre-composition with the functors $\iota_I^n$ and $\pi_I^n$
form a pair of left and right adjoint functors.

By definition of the cross effect of a functor we get
$$G(V_1\oplus\dots\oplus V_n)=\bigoplus_{m=1}^n\bigoplus_{i_1<\dots <i_m} G^{(m-1)}(V_{i_1}|\dots|V_{i_m}).$$
Hence
$$G\circ\Pi_n=\bigoplus_{m=1}^n\bigoplus_{i_1<\dots <i_m}G^{(m-1)}\circ\pi_{\{i_1,\dots,i_m\}}^n.$$
Moreover, using that $\widehat{\iota_I^n}$ is left adjoint to $\widehat{\pi_I^n}$,
in case $m\neq n$ we get
\begin{eqnarray*}
Hom_{n-{\cal A}_{R}}((R\otimes\boxtimes^n),
\lefteqn{G^{(m-1)}\circ\pi_{\{i_1,\dots,i_m\}}^n)=}\\&=&Hom_{n-{\cal A}_{R}}((R\otimes\boxtimes^n)\circ\iota_{\{i_1,\dots,i_m\}}^n,
G^{(m-1)})=0
%\\&=&0
\end{eqnarray*}
since at least one of the factors in the tensor product is $0$.
Therefore
$$
Hom_{{\cal A}_{R}}((R\otimes id^{\otimes n},G)=
Hom_{n-{\cal A}_{R}}((R\otimes\boxtimes^n),G^{(n-1)})
$$
Since $G$ is of degree at most $n$, the functor $G^{(n-1)}$ is
additive in each variable. Hence its value is
given by the $R$-module $G^{(n-1)}(\Z|\dots|\Z)$. However,
$(R\otimes id^{\boxtimes^n})(\Z|\dots|\Z)=R$ and
$Hom_R(R,G^{(n-1)}(\Z|\dots|\Z))=G^{(n-1)}(\Z|\dots|\Z)$.
Hence a natural transformation in
$Hom_{n-{\cal A}_{R}}((R\otimes\boxtimes^n),G^{(n-1)})$
induces an element in $G^{(n-1)}(\Z|\dots|\Z)$. On the other hand,
any element in $G^{(n-1)}(\Z|\dots|\Z)$ induces a natural transformation
in $Hom_{n-{\cal A}_{R}}((R\otimes\boxtimes^n),G^{(n-1)})$.
Therefore $$
Hom_{{\cal A}_{R}}(R\otimes id^{\otimes n},G)=G^{(n-1)}
$$
as claimed and
we obtain the statement. \dickebox

\begin{Prop}\label{recollementA}
$$
{\cal A}^{n-1}_{\F_p}{
\begin{array}c\leftarrow\\ \rightarrow\\ \leftarrow \end{array}}
{\cal A}^{n}_{\F_p}{
\begin{array}c\leftarrow\\ \stackrel{e}{\rightarrow}\\ \leftarrow \end{array}}
{\mathbb F}_p{\Sn_n}-mod
$$
is a recollement diagram with
$$e:=Hom_{{\mathcal A}_{\F_p}}(\F_p\otimes_{\Z}\underbrace{id
\otimes_{\Z} \dots
\otimes_{\Z}id}_{\mbox{\footnotesize $n$ factors}},-):
{\cal A}^{n}_{\F_p}\lra {\mathbb F}_p{\Sn_n}-mod.$$
\end{Prop}

Proof.
We shall show the following auxiliary lemma needed for the
proof of the proposition.

\begin{Lemma}
$\F_p\otimes id^{\otimes n}$ is projective
and injective in ${\mathcal A}_{\F_p}^n$.
\end{Lemma}

Proof.
Taking cross effects is exact. By Lemma~\ref{macdonaldlinearisation}
we see that $\F_p\otimes id^{\otimes n}$ is projective.
Using that the duality $D$ on functors $F:\Z-free\lra \F_p-Mod$
given by
$$(DF)(V)=Hom_{\F_p}(F(Hom_\Z(V,\Z)),\F_p),$$
and observing that $\F_p\otimes id^{\otimes n}$ is self-dual, we obtain that
$\F_p\otimes id^{\otimes n}$ is injective as well. \dickebox

\medskip

We need to prove that
$${\cal A}^{n-1}_{\F_p}=\{F\in{\cal A}^{n}_{\F_p}|\;
Hom_{{\mathcal A}_{\F_p}}(\F_p\otimes id^{\otimes n},F)=0\}.$$
But this is clear by Lemma~\ref{macdonaldlinearisation}

Since $\F_p\otimes id^{\otimes n}$ is projective
there is a right adjoint and a left adjoint to
$Hom_{{\mathcal A}_{\F_p}}(\F_p\otimes id^{\otimes n},-),$
namely the functor
$M\mapsto (\F_p\otimes id^{\otimes n}\otimes M)^{\Sn_n}$
is the right adjoint, and the functor
$M\mapsto (\F_p\otimes id^{\otimes n}\otimes M)_{\Sn_n}$
is the left adjoint.

Moreover, the unit and the counit of the adjunctions induce the identity
on $\F_p\Sn_n$-mod. This can be done literally as in Piriou
\cite[Proposition 2.2.2]{piriouthese}.

This shows Proposition~\ref{recollementA}.
\phantom{x}\hfill\dickebox

\bigskip

As a consequence we show the following lemma.

\begin{Lemma}\label{simplescorrespond}
Any simple object in ${\mathcal A}^n_{\F_p}$ is in the image of
$\circ(\F_p\otimes-):
{\mathcal F}^n_{\F_p}\hookrightarrow {\mathcal A}^n_{\F_p}$
and any simple object of ${\mathcal A}^n_{\F_p}$ gives a
simple object in ${\mathcal F}^n_{\F_p}$ this way.
\end{Lemma}

Proof. We shall use induction on $n$. There is a morphism of
recollement diagrams as follows, where the vertical functors are
fully faithful embeddings of categories by Lemma~\ref{inclusion}.
$$
\begin{array}{ccccc}
{\cal F}^{n-1}_{\F_p}&{
\begin{array}c\leftarrow\\ \rightarrow\\ \leftarrow \end{array}}&
{\cal F}^{n}_{\F_p}&{
\begin{array}c\leftarrow\\ \stackrel{e_{{\mathcal F}^n}}{\rightarrow}\\
\leftarrow \end{array}}&
{\mathbb F}_p{\Sn_n}-mod\\ \\
\dar i_{n-1}&&\dar i_n&&\|\\ \\
{\cal A}^{n-1}_{\F_p}&{
\begin{array}c\leftarrow\\ \rightarrow\\ \leftarrow \end{array}}&
{\cal A}^{n}_{\F_p}&{
\begin{array}c\leftarrow\\ \stackrel{e_{{\mathcal A}^n}}{\rightarrow}\\
\leftarrow \end{array}}&
{\mathbb F}_p{\Sn_n}-mod
\end{array}
$$
Therefore the number of simple objects in ${\cal A}^{n}_{\F_p}$ and in
${\cal F}^{n}_{\F_p}$ coincides.

The statement is clear for $n\leq p-1$ by Lemma \ref{charp}.
Let $n\geq p$ and let $S$ be a simple object in ${\mathcal
F}^n_{\F_p}$. We may suppose, using the induction hypothesis that
$e_{{\mathcal F}^n}(S)\neq 0$.
Then, suppose $X$ is a simple and proper subobject of $i_nS$.
Since $e_{{\mathcal A}^n}$ is exact, $e_{{\mathcal A}^n}(X)$ is a
subobject of $e_{{\mathcal A}^n}(i_nS)=e_{{\mathcal F}^n}(S)$.
Since $e_{{\mathcal F}^n}$ is exact, $e_{{\mathcal F}^n}(S)$ is
simple. So, $e_{{\mathcal A}^n}(X)$ is either $0$ or isomorphic to
$e_{{\mathcal F}^n}(S)$. Since $X$ is a proper subobject of $i_nS$
we see that $e_{{\mathcal A}^n}(X)=0$, and hence $X\in{\mathcal
A}^{n-1}$. Since the simple objects of
${\mathcal F}^{n-1}$ and of ${\mathcal A}^{n-1}$ coincide by the
induction hypothesis, $X$ is a proper non zero subobject
of $S$. This gives the contradiction. \hfill\dickebox

\bigskip

We have seen in Lemma~\ref{charp} and Remark~\ref{nkleinerp} that
$${\mathcal F}^{p-1}_{\F_p}\simeq{\mathcal A}^{p-1}_{\F_p}\simeq
\prod_{n<p}\prod_{\lambda \vdash n}\F_p-mod$$

\section{Short Review of Brauer tree algebras}

\label{Brauertreealgebrasection}

Proposition~\ref{recollementA} shows that the representation theory
of the symmetric group is closely related to ${\mathcal A}^n(p)$
we shall need some information from group representations. We give Benson~\cite[Section 4.18 and Section 6.5]{Benson} as a general reference.

Let $k$ be a field of characteristic $p>0$ and let $G$ be a finite group
so that $k$ is a splitting field for $G$, i.e. the endomorphism ring of
each simple $kG$-module is $k$. Then it is well-known that there are
only a finite number of isomorphism classes of indecomposable $kG$-modules
if and only if the Sylow $p$ subgroups of $G$ are all cyclic.
In particular the
symmetric group ${\mathfrak S}_p$ is a group with every Sylow $p$
subgroup being cyclic.

An indecomposable ring-direct factor of $kG$ is called a $p$-block of $G$,
and the representation theory of a block $B$ is completely understood
in case the Sylow $p$ subgroups of $G$ are cyclic.
The algebras occurring in this case are so-called Brauer tree algebras.

We explain briefly the theory of Brauer tree algebras and refer to
Benson~\cite[Section 4.18 and Section 6.5]{Benson} or to Auslander-Reiten-Smal\o{}~\cite[Section X.3]{ARS} for more details.

A Brauer tree is a finite connected tree $\Gamma$, i.e. a connected
graph without cycles and multiple edges, together with additional
structure as we shall describe now:
Let $\Gamma_0$ be the vertices of $\Gamma$ and let $\Gamma_1$ be
the edges of $\Gamma$. For each $e\in\Gamma_1$ let
$v_1(e)\in\Gamma_0$ and $v_2(e)\in\Gamma_0$
be the vertices adjacent to $e$.
\begin{itemize}
\item
Then for each vertex $v\in\Gamma_0$ let
$E_v$ be the edges adjacent to $v$, and suppose given
a transitive permutation $\sigma_v\in{\mathfrak S}_{E_v}$.
\item
Moreover, we choose one vertex $v_0\in\Gamma_v$ and an integer $\mu>0$, called the exceptional
multiplicity.
\end{itemize}
A Brauer tree algebra $B$ is then a finite dimensional symmetric
$k$-algebra so that
\begin{itemize}
\item
the isomorphism classes of simple $B$-modules
are parameterised by $\Gamma_1$, denoting by $S_e$ a simple $B$-module
whose isomorphism class corresponds to the edge $e\in\Gamma_1$,
\item
the projective cover $P_e$ of $S_e$ has the property that $rad(P_e)/soc(P_e)=U_{v_1(e)}\oplus U_{v_2(e)}$, where $U_{v_1(e)}$
and $U_{v_2(e)}$ are both uniserial modules.
\item
for each $j\in\{1,2\}$,
\begin{itemize}
\item if $v_j(e)$ is not the exceptional vertex, then
the composition length of $U_{v_j(e)}$ is equal to $|E_{v_j(e)}|-1$
\item
if $v_j(e)$ is the exceptional vertex, then
the composition length of $U_{v_j(e)}$ is equal to $\mu\cdot |E_{v_j(e)}|-1$
\end{itemize}
\item
for each $i$, for which the $i$-th radical quotient is not $0$,
$$rad^i(U_{v_1(e)})/rad^{i+1}(U_{v_1(e)})\simeq S_{\sigma_{v_1(e)}^{i+1}(e)}$$
and
$$rad^i(U_{v_2(e)})/rad^{i+1}(U_{v_2(e)})\simeq S_{\sigma_{v_2(e)}^{i+1}(e)}$$
\end{itemize}

One observes that in case $\mu=1$, then there is no difference between
the exceptional vertex and a non exceptional vertex. Hence, in this case we
do not need to fix an exceptional vertex and we say that the Brauer tree has
no exceptional vertex if the exceptional multiplicity $\mu$ is $1$.

Observe further that if a vertex $v$ is a leaf, i.e. $E_v=\{e\}$, then
the projective cover of $S_e$ is uniserial and $U_v=0$.

A particular case is when each vertex has at most $2$ edges adjacent to it,
or in other words if $|E_v|\leq 2$ for all $v\in\Gamma_0$.
In this case we call the Brauer tree a stem. If moreover the Brauer tree
 has no exceptional vertex then the tree can be visualised as
$$\bullet-\bullet-\bullet-\dots-\bullet-\bullet$$
The projective indecomposable modules $P$ are then all of Loewy length $3$,
i.e. $P$ has simple top, simple socle and $rad(P)/soc(P)$ is semisimple of composition length at most $2$.

This situation occurs for $\F_p{\mathfrak S}_p$ which has several blocks for $p\geq 5$,
one of which is a Brauer tree algebra associated to a stem
without exceptional vertex and $p-1$ edges,
and all the other blocks are simple algebras.

\section{The structure of polynomial functors modulo $p$ }
\label{modpcase}
The situation of polynomial functors of degree $p$ is different
from those of degree $n<p$. We are going to describe in this section
their structure completely. From now on we assume that $p\geq 5$
since the representation theory of $\F_2\Sn_2$ and of $\F_3\Sn_3$
is slightly different from the case $p\geq 5$.

\label{SchurBrauer}

\begin{Rem}\label{structureofschuralgebra}
Let us recall the relation between the Schur
algebra $S_{\F_p}(p,p)$ and the group algebra $\F_p\Sn_p$. We shall use
the structure of the
Schur algebra and the corresponding Hecke algebra
as it is proved in \cite[Corollary 1.3]{cyclotomic}. The
group algebra of the symmetric group is a special case of the Hecke algebra.
For general
informations on Schur algebras and Brauer tree algebras see
\cite{Green} and \cite{derbuch}.
The
algebra $S_{\F_p}(p,p)$ admits $p$ simple modules
$S_1,\dots,S_{p-1},S_p$
whereas the group algebra $\F_p\Sn_p$ admits $p-1$ simple modules
$S_1',\dots,S_{p-1}'$. The projective indecomposable
$S_{\F_p}(p,p)$-modules have composition series
$$
\begin{array}{c}
S_1\\S_2\\S_1
\end{array},\;
\begin{array}{ccc}
&S_2\\S_1&&S_3\\&S_2
\end{array},\;
\begin{array}{ccc}
&S_3\\S_2&&S_4\\&S_3
\end{array},
\;\dots,\;
\begin{array}{ccc}
&S_{p-1}\\S_{p-2}&&S_p\\&S_{p-1}
\end{array},\;
\begin{array}{c}
S_p\\S_{p-1}
\end{array}
$$
whereas the projective indecomposable $\F_p\Sn_p$-modules
have composition series
$$
\begin{array}{c}
S_1'\\S_2'\\S_1'
\end{array},\;
\begin{array}{ccc}
&S_2'\\S_1'&&S_3'\\&S_2'
\end{array},\;
\begin{array}{ccc}
&S_3'\\S_2'&&S_4'\\&S_3'
\end{array},\;\dots,\;
\begin{array}{ccc}
&S_{p-2}'\\S_{p-3}'&&S_{p-1}'\\&S_{p-2}'
\end{array},\;
\begin{array}{c}
S_{p-1}'\\S_{p-2}'\\S_{p-1}'
\end{array}
$$
\end{Rem}

\begin{Lemma}\label{identifyingtheext}
Let $L$ be the simple polynomial functor in ${\mathcal F}_{\F_p}^p$
so that $L$ is of degree $p$ and so that $L$
corresponds to the trivial representation of $\F_p\Sn_p$. Then,
$Ext^1_{{\mathcal F}_{\F_p}}(id,L)\neq 0\neq Ext^1_{{\mathcal
F}_{\F_p}}(L,id)$. Moreover, if
$Ext^1_{{\mathcal F}_{\F_p}}(id,S)\neq 0$ or $0\neq Ext^1_{{\mathcal
F}_{\F_p}}(S,id)$ for a simple degree $p$-functor $S$, then $L\simeq S$.
\end{Lemma}

Proof. Let $S_1$ be the
simple polynomial functor of degree $1$. The
identity functor $id$ is trivially of degree $1$ and simple,
which implies $S_1=id$. Then, $Ext_{{\cal F}^p}^1(L,id)$ is not
necessarily zero for $L$ being an irreducible polynomial functor
of degree $p$. Now, since we are working over $\F_p$, we get that
$id^{(1)}\simeq id$ as polynomial functor, but not as strict
polynomial functor. As strict polynomial functor, $I^{(1)}$ is of
degree $p$.

Theorem \ref{Kuhnresult} implies that there is only
one simple functor $L$ of degree $p$ with
$Ext_{{\cal F}^p}^1(L,id)\neq 0$.
Proposition \ref{Andersen} in connection with
Theorem \ref{FFSSresult}, imply that $L$ is the simple functor
corresponding to the trivial $\F_p\Sn_p$-module, since this is the
module which has an extension with the unique simple module of the
Schur algebra $S_{\F_p}(p,p)$ which is not a simple
$\F_p\Sn_p$-module (cf Remark~\ref{structureofschuralgebra}).

This
implies that
$$Ext_{{\cal P}(p)}^1(L,id^{(1)}) \simeq
Ext_{{\cal F}_{\F_p}}^1(L,id^{(1)})\simeq
Ext_{{\cal F}_{\F_p}}^1(L,id)\;.$$
Since the
category of degree $p$ strict polynomial functors is equivalent to
the category of modules over the Schur algebra $S_{\F_p}(p,p)$, one sees
that
$$Ext_{{\cal P}_{\F_p}}^1(L,id^{(1)})\simeq Ext_{S_{\F_p}(p,p)}^1(V,I_0)$$
for $I_0$ being the simple $S_{\F_p}(p,p)$-module
corresponding to the $p$ singular partition of $p$ and $V$ being
the simple $S_{\F_p}(p,p)$-module corresponding to $L$. Finally, it is a
classical fact (cf e.g. \cite{cyclotomic})
that $Ext_{S_{\F_p}(p,p)}^1(I_0,V)\neq 0$ or
$Ext_{S_{\F_p}(p,p)}^1(V,I_0)\neq 0$ if and only if $V$ corresponds to
the trivial module of the symmetric group, and in this case, the
dimension of $Ext_{S_{\F_p}(p,p)}^1(I_0,V)$ and  of
$Ext_{S_{\F_p}(p,p)}^1(V,I_0)$ is $1$. This proves the statement.
\phantom{x}\hfill\dickebox

\medskip

Using the embedding ${\mathcal F}^p_{\F_p}\hookrightarrow
{\mathcal A}^p_{\F_p}$ from Lemma~\ref{inclusion},
the functor $L$ of ${\mathcal F}^p_{\F_p}$ induces a functor
$L(\F_p\otimes-)$ in ${\mathcal A}^p_{\F_p}$. In order to avoid
additional notational burden we shall denote this functor $L(\F_p\otimes-)$
by $L$ as well.

We get as a corollary the following statement.

\begin{Cor}\label{identifyingL}
$Ext^1_{{\mathcal A}^p_{\F_p}}(\F_p\otimes_\Z id,L)\neq 0\neq
Ext^1_{{\mathcal A}^p_{\F_p}}(L,\F_p\otimes_\Z id)$.
\end{Cor}

Proof.
We know that $Ext^1_{{\mathcal F}^p_{\F_p}}(id,L)\neq 0$.
So, there is a non split exact sequence
$$0\lra L\lra X\lra id\lra 0$$ for some functor
$X$ in ${\mathcal F}^p_{\F_p}$.
Since ${\mathcal F}^p_{\F_p}\hookrightarrow
{\mathcal A}^p_{\F_p}$ by Lemma~\ref{inclusion}.
This induces an exact sequence
$$0\lra L(\F_p\otimes_{\Z}id)\lra
X(\F_p\otimes_\Z id)\lra \F_p\otimes_{\Z}id\lra 0$$
in ${\mathcal A}^p_{\F_p}$ where $L(\F_p\otimes_{\Z}id)$ is simple by
Lemma~\ref{simplescorrespond}. This sequence is non split
since the functor pre-composing with $\F_p\otimes_\Z id$
is a fully faithful embedding.
Hence $Ext^1_{{\mathcal A}^p_{\F_p}}(\F_p\otimes_\Z id,L)\neq 0$.
Therefore, $L$
is a direct factor of the top of the radical of the projective cover of
$\F_p\otimes_{\Z}id$. Similarly,
$Ext^1_{{\mathcal A}^p_{\F_p}}(L,\F_p\otimes_\Z id)\neq 0$.
\hfill\dickebox

\medskip

Actually, the argument of Lemma~\ref{identifyingtheext}
gives another slightly different statement.

\begin{Lemma} Let $X$ and $Y$ be two simple functors of degree at most
$p$. Then, $$Ext^1_{{\mathcal F}_{\F_p}^p}(X,Y)\neq 0\Rightarrow
deg(X)-deg(Y)\in\{0,p-1\}.$$
\end{Lemma}

Proof. By Lemma~\ref{simplescorrespond} we know that the simple functors
$X$ and $Y$ can be considered to lie in ${\mathcal F}_{\F_p}^p$.
Proposition \ref{Andersen} in connection with
Theorem \ref{FFSSresult} imply this result.
\hfill\dickebox

\begin{Rem}
At the present stage
it might happen that $Ext^1_{{\mathcal A}_{\F_p}^p}(X,Y)\neq 0$
even though $Ext^1_{{\mathcal F}_{\F_p}^p}(X,Y)= 0$.
\end{Rem}

Denoting by $\rho_p(i)$
the number of $p$-regular partitions of $i$, the algebra
$\Gamma_{\F_p}^p$ is Morita equivalent to a direct product of
$\left(\sum_{i=1}^{p-1}\rho_p(i)\right)-1$ copies of $\F_p$ and of an
indecomposable ring
$\Gamma_{\F_p,0}^p$. By the recollement diagram preceding
Lemma~\ref{simples} this ring
$\Gamma_{\F_p,0}^p$ has a projective module
$P=\Gamma_{\F_p,0}^p\cdot e$ so that the endomorphism ring of $P$
is  Morita equivalent to the Brauer tree algebra corresponding to
$\F_p\Sn_p$.
There is a projective indecomposable $\Gamma_{\F_p,0}^p$-module $P_0$
so that $P_0\oplus P$ is a progenerator of $\Gamma_{\F_p,0}^p$
and the endomorphism ring of $P_0\oplus P$ is basic and Morita
equivalent to $\Gamma_{\F_p,0}^p$. Moreover,
$\Gamma_{\F_p,0}^p/(\Gamma_{\F_p,0}^p\cdot e\cdot \Gamma_{\F_p,0}^p)$
is Morita equivalent to $\F_p$.

\medskip

The next remark constructs appropriate projective objects.

Let $R$ and $S$ be commutative rings.
The functor $R[-]/(I^{n+1})$ which assigns to an $S$-module
$V$ the quotient of the semi-group ring $R[V]$ on $V$ by
the $n+1$-st power of the augmentation ideal. This functor
$S-free\lra R-mod$ is polynomial of degree $n$. Define
$proj_{n}^{m}:=\F_p[Hom_\Z(\Z^m,-)]/(I^{n+1})$
and $proj_\infty^m:=\F_p[Hom_\Z(\Z^m,-)]$.

\begin{Lemma}
The functor $proj_n^m$ in ${\mathcal A}_{\F_p}^n$ is projective
and contains a projective cover of the reduction modulo $p$ functor
$\F_p\otimes_{\Z}id$.
\end{Lemma}

Proof. For any degree $n$ polynomial functor $F$ one gets
$Hom_{{\mathcal A}_{\F_p}}(proj_n^m,F)\simeq F(\Z^m)$.
This is true if one does
not factorizes the power of the augmentation ideal, and since all
functors are of degree at most $n$, each natural transformation
from $proj_{\infty}^m$ to $F$ is zero on $I^{n+1}$ (see
\cite[Section 1]{FLS}).

So, $Hom_{\mathcal A}(proj_n^m,-)$ is exact,
as evaluation on exact sequences of functors
is exact. Hence, $proj_n^m$ is a projective object in
${\mathcal A}_{\F_p}^n$ and since
$Hom_{\mathcal A}(proj_n^m,id)=id(\Z^m)=\Z^m\neq 0$, the
projective cover of the reduction modulo $p$
functor is a direct summand in $proj_n^m$.\hfill\dickebox

\begin{Rem}
The situation is different for ${\mathcal F}_{\F_p}$. Indeed,
the  functor $\F_p[Hom_{\F_p}(\F_p^m,-)]/I^{n+1}$ is projective
in ${\mathcal F}_{\F_p}^n$ again since
$$Hom_{{\mathcal F}_{\F_p}^n}(\F_p[Hom_{\F_p}(\F_p^m,-)]/I^{n+1},F)
\simeq F(\F_p^m).$$ But, automatically $I^{p}=0$ for $m=1$ and
the evaluation at $\F_p$ in this case, so that
the endomorphism ring of $\F_p[Hom_{\F_p}(\F_p,-)]/I^{n+1}$
is an $\ell$-dimensional vector space where $\ell=min(p,n)$.

Recall the embedding ${\mathcal F}_{\F_p}^n\lra {\mathcal A}_{\F_p}^n$
given by pre-composing with $\F\otimes_\Z id$. The image of
$\F_p[Hom_{\F_p}(\F_p,-)]/I^{n+1}$ under this embedding is
$\F_p[Hom_{\F_p}(\F_p,\F_p\otimes_\Z -)]/I^{n+1}$ which is
different from $\F_p[Hom_{\Z}(\Z,-)]/I^{n+1}$. As we will see,
the projective indecomposable cover of $\F_p\otimes_\Z id$
in ${\mathcal A}_{\F_p}^p$ is a direct factor
of the functor $\F_p[Hom_{\Z}(\Z,-)]/I^{p+1}$
and this is the only indecomposable functor which is not in
the image of the embedding
${\mathcal F}_{\F_p}^p\lra {\mathcal A}_{\F_p}^p$.
\end{Rem}

\begin{Lemma}\label{compofatorsofproj}
The projective functor $proj_n^1$ has one composition factor of
degree $d$ for each $0\leq d\leq n$ for all $n\leq p-1$.
In particular $proj^1_n$ contains a simple constant functor as
a direct summand.
\end{Lemma}

Proof. First, $proj^1_n(0)=\F_p$, and so the simple functor
of degree $0$ is a direct factor of $proj^1_n$.

Furthermore, $proj_n^1\epi proj_{n-1}^1$ for trivial reasons.
Moreover, $proj_n^1$ is of degree $n$ and not of degree $n-1$.
We compute that $End_{{\mathcal
A}_{\F_p}^n}(proj_n^1)=proj_n^1(\Z)=\F_p[\Z]/I^{n+1}$
is an $n+1$-dimensional vector space.
Moreover, by Lemma~\ref{charp} and Remark~\ref{nkleinerp}
the projective module $proj_{p-1}^1$ is semisimple since
${\mathcal A}_{\F_p}^{p-1}$ is a semisimple category.
Observe that $proj_n^1$ has
exactly one composition factor more than $proj_{n-1}^{1}$. This
composition factor is of degree $n$ since
$$Hom_{{\mathcal A}_{\F_p}^n}(proj_n^1,\F_p\otimes_{\Z}\underbrace{id
\otimes_{\Z}\dots
\otimes_{\Z}id}_{\mbox{\footnotesize $n$ factors}})=
\F_p\otimes_{\Z}\Z\otimes_{\Z}\dots\otimes_{\Z}\Z=\F_p$$
and since $\F_p\otimes_{\Z}\underbrace{id
\otimes_{\Z}\dots
\otimes_{\Z}id}_{\mbox{\footnotesize $n$ factors}}$ is the
projective object corresponding to degree $n$ polynomial functors in the
recollement diagram.
\hfill\dickebox

\medskip

Recall that $L$ denotes the simple functor in ${\mathcal A}^p_{\F_p}$
mapping to the trivial $\F_p\Sn_p$-module in the recollement
diagram.

\begin{Prop}\label{structureofP0} Suppose $p\geq 5$.
The projective cover $P_{\F_p\otimes id}$ of $\F_p\otimes_{\Z}id$
in ${\mathcal A}_{\F_p}^p$
is uniserial with top and socle being $\F_p\otimes_{\Z}id$ and with
$rad(P_{\F_p\otimes_\Z id})/soc(P_{\F_p\otimes_\Z id})\simeq L$.
Moreover, $$proj^1_p\simeq S_0\oplus P_{\F_p\otimes id}\oplus
S_2\oplus\dots \oplus S_{p-1}$$ for simple functors $S_i$ of degree $i$.
\end{Prop}

Proof.
We shall divide the proof into various claims.

\begin{Claim}
No direct summand of the top of $proj^1_p$ is of degree $p$.
\end{Claim}

Proof: We have $Hom_{{\mathcal
A}_{\F_p}^p}(proj_p^1, (\F_p\otimes_\Z id)^{\otimes p})=\F_p$.
On the other hand, we know from the recollement diagram
that $(\F_p\otimes_\Z id)^{\otimes p}$ is the
projective cover of the simple modules coming from
$\F_p\Sn_p$, that is those of degree $p$.
Now, each projective indecomposable module of $\F_p\Sn_p$ has the
property that the top of this module is isomorphic to the socle of
this module and that the top and the socle of this projective
indecomposable module are different. Hence, suppose a simple polynomial
functor of degree $p$ would be in the top of $proj^1_p$, then
let $Q$ be its projective cover in ${\mathcal A}_{\F_p}^p$. Further,
$Q$ is a direct summand of $(\F_p\otimes_\Z id)^{\otimes p}$.
Since the top of $(\F_p\otimes_\Z id)^{\otimes p}$ is isomorphic to
the socle of $(\F_p\otimes_\Z id)^{\otimes p}$,
the above homomorphism space would be at least $2$-dimensional,
corresponding to the mapping of $proj^1_p$ on the top and on 
the socle of  $(\F_p\otimes_\Z id)^{\otimes p}$.
\hfill\dickebox

\medskip

\begin{Lemma}\label{degreesunderextinAp}
If $Ext^1_{{\mathcal A}_{\F_p}^p}(S,T)\neq 0$ for two simple functors
$S$ and $T$, then $deg(S)-deg(T)\in\{0,p-1\}$ and if $deg(S)=deg(T)$, then
$deg(S)=p$.
\end{Lemma}

Proof.
We know by Lemma~\ref{charp} that ${\mathcal A}_{\F_p}^{p-1}$
is semisimple. Moreover, the category of constant functors is
a direct factor in the category of polynomial functors. Using
Lemma~\ref{degreeunderext} this shows the statement.\hfill\dickebox

\medskip

We denote by $P_V$ the projective cover of the functor $V$ in
${\mathcal A}^p_{\F_p}$.

\begin{Claim}\label{structureofM}
$proj_p^1\simeq S_0\oplus S_2\oplus S_3\oplus
\dots\oplus S_{p-1}\oplus M$
for simple projective functors $S_i$ of degree $i$ and the projective
cover $M$ of the functor $\F_p\otimes_{\Z}id$.
\end{Claim}

Proof.
$End_{{\mathcal A}_{\F_p}^p}(proj_p^1)=proj_p^1(\Z)=\F_p[\Z]/I^{p+1}$
is a $p+1$-dimensional $\F_p$-vector space.
We know already that $proj_{p-1}^1$ is a quotient of
$proj_p^1$ and that this is a semisimple functor with $p-1$
direct factors. So, every direct summand of the semisimple functor
$proj^1_{p-1}$ is a direct factor of the head of $proj^1_p$.
Denote by $S_0,S_1,S_2\cdots S_{p-1}$ the simple direct factors of
$proj^1_{p-1}$ and let $S_i$ be of degree $i$. Then, since the degree $0$
functors split off in any case, $P_{S_0}=S_0$. Moreover,
$S_0\oplus P_{S_0}\oplus P_{S_1}\oplus \dots\oplus P_{S_{p-1}}$ is a
direct factor of $proj^1_p$.

We need to study the functors $P_{S_i}$. If all the composition
factors of $P_{S_i}$ for an $i\leq p-1$ are of degree $p-1$ at most, then
by Lemma~\ref{degreeunderext} we get that $P_{S_i}$ is of degree at most
$p-1$ as well. Since the category ${\mathcal A}_{\F_p}^{p-1}$ is semisimple,
we get $P_{S_i}=S_i$.

Suppose that a degree $p$ simple polynomial functor $S$ is a composition
factor of $P_{S_{i_0}}$. Then, since ${\mathcal A}_{\F_p}^{p-1}$ is
semisimple, again by Lemma~\ref{degreeunderext},
$Ext^1_{{\mathcal A}_{\F_p}^{p}}(S_{i_0},S)\neq 0$.
Now, simple functors are self-dual under the duality
$(DF)(V):=F(V^*)^*$ (cf \cite{KuhnII} for functors in ${\mathcal F}_{\F_p}$
and by Lemma~\ref{simplescorrespond} for simple functors in
${\mathcal A}_{\F_p}$).
So, $Ext^1_{{\mathcal A}_{\F_p}^{p}}(S,S_{i_0})\neq 0$.
Since $Hom_{{\mathcal A}_{\F_p}^{p}}(proj^1_p,\F_p\otimes_{\Z}-)$ is
a one-dimensional $\F_p$-vector space,
this happens for precisely one $i_0\in \{1,2,\dots,p-1\}$. We already
know by Claim~\ref{identifyingL}
that $i_0=1$ and $S$ is the simple module corresponding to the trivial
representation of $\F_p\Sn_p$. This proves the claim. \hfill\dickebox

\begin{Claim}\label{structureofPL}
For the projective cover $P_L$ of $L$ we get
$rad(P_L)/rad^2(P_L)\simeq (\F_p\otimes id)\oplus L_2$, where $L_2$ is simple of
degree $p$, $L_2\not\simeq L$ and $rad^2(P_L)\simeq L$ as well as $rad^3(P_L)=0$.
\end{Claim}

\begin{Rem}\label{RemarkstructureofPl}
We do not claim here that $soc(P_L)$ is simple. However, the radical layer
structure of $P_L$ can be described by
$$\begin{array}{ccc}&L\\ (\F_p\otimes id)&&L_2\\&L\end{array}$$
and where it is not clear if $soc(P_L)$ is simple and isomorphic to $L$
or if the socle is isomorphic to $L\oplus L_2$ or to $L\oplus  \F_p\otimes id$.
But
$Hom_{{\mathcal A}_{\F_p}^p}(\F_p\otimes id^{\otimes p},P_L)$
is the projective cover of the trivial $\F_p\Sn_p$-module. This
projective cover is uniserial with composition series
$$\begin{array}{c}L\\ L_2\\ L\end{array}.$$
Since $L_2$ is simple of degree $p$, its image in $\F_p\Sn_p$ is given by
the known module structure of $\F_p\Sn_p$. In particular, $L\not\simeq L_2$.
This shows that the only uniserial module of length $3$
which is a quotient of $P_L$, if there is any,
can have composition series $$\begin{array}{c}L\\ L_2\\ L\end{array}.$$
In particular, $soc(P_L)\not\simeq L\oplus L_2$.
\end{Rem}

Proof of Claim~\ref{structureofPL}. By Claim~\ref{identifyingL} we know that
$(\F_p\otimes id)$ is composition factor of $top(rad(P_L))$ and since
$Hom_{{\mathcal A}_{\F_p}^{p}}(proj^1_p,\F_p\otimes_{\Z}-)$ is
a one-dimensional $\F_p$-vector space we know that it has
multiplicity $1$.
Since the image of $P_L$ in $\F_p\Sn_p-mod$ is uniserial with top and
socle $L$ and simple $rad(P_L)/soc(P_L)\simeq L_2$, we have the above
structure.\hfill\dickebox

\begin{Claim}\label{topradM}
The projective cover $M=P_{\F_p\otimes id}$
of $\F_p\otimes_\Z id$ has $rad(M)/rad^2(M)\simeq L$.
\end{Claim}

Proof. By Claim~\ref{identifyingL} we know that $L$ is a direct factor
of $top(rad(M))$. Since ${\mathcal A}_{\F_p}^{p-1}$ is semisimple,
using Lemma~\ref{degreeunderext} we see that no simple functor of degree
$p-1$ at most can be a direct factor of $top(rad(M))$.
Suppose $top(rad(M))$ has a second simple direct factor $T$ of degree $p$.
Then, $Ext^1_{{\mathcal A}_{\F_p}^{p}}(\F_p\otimes_\Z id, T)\neq 0$.
Simple functors are self-dual (cf as above \cite{KuhnII} for
functors in ${\mathcal F}_{\F_p}$
and by Lemma~\ref{simplescorrespond} for simple functors in
${\mathcal A}_{\F_p}$). So,
$Ext^1_{{\mathcal A}_{\F_p}^{p}}(T,\F_p\otimes_\Z id)\neq 0$.
But, we have seen in Corollary~\ref{identifyingL}
that there is one simple functor of degree $p$
with a non trivial extension group with $\F_p\otimes_\Z id$, namely $L$.
Moreover, since $proj^1_p$ contains the projective cover of the
simple degree $1$-functor as a direct factor
(see Claim~\ref{structureofM}), and since by
Lemma~\ref{degreesunderextinAp} this is the only degree where
non trivial first extension groups can occur, we see that $T\simeq L$.

Suppose $ L\oplus L| rad(M)$.
Then, $Ext^1_{{\mathcal A}_{\F_p}^{p}}(\F_p\otimes_\Z id, L)$
is two-dimensional at least,
and again by the self-duality of the simple functors,
$Ext^1_{{\mathcal A}_{\F_p}^{p}}(L,\F_p\otimes_\Z id)$ is
at least two-dimensional.
Therefore, $\F_p\otimes_\Z id$ occurs twice in $top(rad(P_L))$.
Since $M$ is a direct factor of $proj^1_p$, the space
$Hom_{{\mathcal A}_{\F_p}^{p}}(proj^1_p,(\F_p\otimes_\Z id)^{\otimes^p})$
would be two-dimensional at least. This contradiction shows that $rad(M)$
has simple top $L$. Hence, $top(rad(M))\simeq L$.
\hfill\dickebox

\begin{Rem}
The radical layer structure of $M$ is therefore given by
$$\begin{array}{c} (\F_p\otimes id)\\ L\\ rad^2(M)\end{array}$$
and by Claim~\ref{structureofPL} we get $rad^3(P_L)=0$, and
therefore we obtain $rad^4(M)=0$.
\end{Rem}

\begin{Claim}\label{radsquaredM}
For the projective cover $M$ of $\F_p\otimes_\Z id$ we get that
$top(rad^2(M))\simeq \F_p\otimes_\Z id$.
\end{Claim}

Proof.
Degree $p$ functors can only have extensions with
degree $p$-functors or degree $1$-functors by
Lemma~\ref{degreesunderextinAp}.
Moreover, the structure of $P_L$ implies that
we get that $top(rad^2(M))$ is a direct summand of $top(rad(P_L))$, whence
is isomorphic to either $0$, or to $L_2$ (which is defined
in Claim~\ref{structureofPL}),
or to $\F_p\otimes_\Z id$, or to $L_2\oplus (\F_p\otimes_\Z id)$.

\medskip

{\bf Suppose $L_2\oplus (\F_p\otimes_\Z id)\simeq top(rad^2(M))$.}
We shall use the fact that by Claim~\ref{structureofPL}
we know the structure of $P_L$.

We get two possibilities for the projective resolution of $\F_p\otimes_\Z id$. 
Either
$$P_L\hookrightarrow M\epi(\F_p\otimes_\Z id)$$
is exact, or
$$L\hookrightarrow P_L\lra  M\epi(\F_p\otimes_\Z id)$$
is exact.

In the second case,
$Ext^2_{{\mathcal A}_{\F_p}^p}((\F_p\otimes_\Z id),L)\neq 0$. By
the self-duality of the simple functors,
$Ext^2_{{\mathcal A}_{\F_p}^p}(L,(\F_p\otimes_\Z id))\neq 0$.
Our information is sufficient for being able to write down
the first terms of the projective resolution of $L$,
$$0\lla L\lla P_L\lla M\oplus P_{L_2}\lla P_L\oplus P_{L_3}\lla\dots$$
for some projective $P_{L_3}$, for some simple object $L_3$ of degree $p$,
given by the known projective resolution
of the trivial $\F_p\Sn_p$-module. Since $p\neq 2$, we get $L_3\not\simeq L$.
In an case $Hom(P_L\oplus P_{L_3},(\F_p\otimes_\Z id))= 0$, and therefore
$Ext^2_{{\mathcal A}_{\F_p}^p}((\F_p\otimes_\Z id),L)= 0$.
This contradiction excludes the case
$$L\hookrightarrow P_L\lra  M\epi(\F_p\otimes_\Z id)$$
is exact.

If $P_L\hookrightarrow M\epi(\F_p\otimes_\Z id)$
is exact, the projective dimension of $(\F_p\otimes_\Z id)$ is $1$.
But, we know by Theorem~\ref{FLS} and the example following it, that
$Ext^2_{{\mathcal F}_{\F_p}^p}(id,id)\neq 0$. By consequence, also
$Ext^2_{{\mathcal A}_{\F_p}^p}
((\F_p\otimes_\Z id),(\F_p\otimes_\Z id))\neq 0$
and therefore the projective dimension of $(\F_p\otimes_\Z id)$ is at
least $2$.

These two observations exclude
$L_2\oplus (\F_p\otimes_\Z id)\simeq top(rad^2(M))$.

\medskip

{\bf Suppose $L_2 \simeq top(rad^2(M))$.}
Then, using the structure of
$P_L$, we get either $rad^3(M)=0$ or $rad^3(M)=L$, and then $rad^4(M)=0$.

If $rad^3(M)=0$, then we get a non split exact sequence
$$0\lla (\F_p\otimes_\Z id)\lla M\lla P_L\lla M\lla L_2\lla 0$$
and therefore
$Ext^3_{{\mathcal F}_{\F_p}^p}(L_2,(\F_p\otimes_\Z id))\neq 0$.
Dualizing,
$Ext^3_{{\mathcal F}_{\F_p}^p}( (\F_p\otimes_\Z id),L_2)\neq 0$.
Our knowledge of the various projective covers of simples is sufficient to
write down the first terms of the projective resolution of
$L_2$. We get
$$0\lla L_2\lla P_{L_2}\lla P_L\oplus P_{L_3}\lla
M\oplus P_{L_2}\oplus P_{L_4}\lla P_L\oplus P_{L_3}\oplus P_{L_5}\lla\dots$$
for projective objects $P_{L_4}$ and $P_{L_5}$
corresponding to degree $p$ simple functors $L_4$ and $L_5$, given by
the known projective resolution of the trivial $\F_p\Sn_p$-module.
Moreover, $L_3$ and $L_5$ are both different from $L_2$, since $p\geq 5$.
This implies
$Ext^3_{{\mathcal F}_{\F_p}^p}( (\F_p\otimes_\Z id),L_2)= 0$.
This contradiction excludes this case as well.

So, assume $rad^3(M)=L$ and $rad^4(M)=0$. This is impossible since
then the endomorphism ring of $M$ would be one-dimensional. This
contradicts the fact that $End_{{\mathcal A}_{\F_p}^p}(proj^1_p)$
is $p+1$-dimensional.

Hence, $L_2\not\simeq top(rad^2(M))$.

\medskip

{\bf We still have the possibility that $rad^2(M)=0$.} But again, this would
imply that $End_{{\mathcal A}_{\F_p}^p}(M)$ would be one-dimensional and
therefore $End_{{\mathcal A}_{\F_p}^p}(proj^1_p)$
is $p$-dimensional. Contradiction.

This proves the claim.\hfill\dickebox

\begin{Claim}
$rad^3(M)=0$.
\end{Claim}

Proof. We know by Claim~\ref{structureofPL} that $rad^3(P_L)=0$.
Since by Claim~\ref{topradM} we have $top(rad(M))\simeq L$,
one sees that $rad^4(M)=0$. Moreover, $rad^3(M)$ is either $0$ or
$L$, since $top(rad^2(M))\simeq (\F_p\otimes_\Z id)$ by
Claim~\ref{radsquaredM} and $top(rad(M))\simeq L$ by Claim~\ref{topradM}.

Suppose $rad^3(M)\simeq L$. Then $M$ is uniserial with
composition length $4$, and $top(rad(M))=L$.
Therefore $P_L$ maps onto $rad(M)$ with image being a uniserial module
$N$ of length $3$ with $rad(N)/soc(N)=\F_p\otimes id$.

But this contradicts the structure of
$P_L$ as described in Claim~\ref{structureofPL} and in particular
Remark~\ref{RemarkstructureofPl}.

This proves the claim.\hfill\dickebox

\medskip

Examining what we showed implies that $proj^1_p$ is as stated
in  Proposition~\ref{structureofP0}. \hfill\dickebox

\bigskip

We now come to our first main result
in describing  the structure of $\Gamma_{\F_p,0}$.
For the relevant definitions on Brauer tree algebras we refer to Section~\ref{Brauertreealgebrasection}.

\begin{Theorem}\label{gammapstructure}
$\Gamma_{\F_p,0}$ is a Brauer tree algebra over $\F_p$ without
exceptional vertex and associated to a stem with $p$ edges.
$$\bullet_1-\bullet_2-\bullet_3-\dots-\bullet_{p+1}$$
\end{Theorem}

Proof. The case $p\leq 3$ is a consequence of
Drozd's results.
By Proposition~\ref{structureofP0} we know that the projective
cover of the functor $\F_p\otimes_{\Z}-$ is uniserial with top and socle
$\F_p\otimes_{\Z}-$ and with second layer $V$, where $V$ is the
simple functor corresponding to the trivial $\F_p\Sn_p$-module.

We know furthermore that except the projective cover of
$\F_p\otimes_{\Z}-$ only the projective indecomposable functor
$P_L$ has a composition factor $\F_p\otimes_{\Z}-$ and that this
composition factor is a direct summand of $top(rad(P_L))$.

Since we know that the principal block of $\F_p\Sn_p$ is a Brauer tree
algebra without exceptional vertex associated to a stem with $p$
vertices,
this means that the we only need to show that $\F_p\otimes_{\Z}-$
is not in the socle of $P_L$, since the only basic algebra with the
composition series as a Brauer tree algebra associated to a stem is 
actually a Brauer tree algebra associated to a stem.

For this we use the duality $D$ on the category of polynomial functors.
The projective indecomposable functor $P_L$ is a direct factor of
$\F_p\otimes_{\Z}
\underbrace{id\otimes_{\Z}\dots
\otimes_{\Z}id}_{\mbox{\small $p$ factors}}$, since this is the
projective cover of all the degree $p$ simple functors.

It is clear that $\F_p\otimes_{\Z}
\underbrace{id\otimes_{\Z}\dots
\otimes_{\Z}id}_{\mbox{\small $p$ factors}}$ is self dual. Since all the
simple functors are self-dual, also $DP_L\simeq P_L$. If
$\F_p\otimes_{\Z}-$ is in the socle of $P_L$, the simple functor
$D(\F_p\otimes_{\Z}-)\simeq \F_p\otimes_{\Z}-$ is in the top of
$DP_L\simeq P_L$, but the top of $P_L$ is $L$ by definition.

This proves the Theorem.\hfill\dickebox

\section{Lifting to characteristic $0$}

\label{lifttochar0}

\subsection{Lifting Brauer tree algebras to orders}

\label{Brauertrees}

\begin{Prop} \label{BBalg}
Let $R$ be a complete discrete valuation ring with
residue field $k$ and field of fractions $K$.
Let $B$ be a Brauer tree algebra over $k$ associated to a
Brauer tree which is a stem without exceptional vertex.
Let $\Lambda$ be an 
$R$-order.
Then, for any proper two-sided ideal $I\neq 0$ of $B$
we get that
$$\Lambda\otimes_Rk\simeq B/I\Longrightarrow
rank_{\Z}(K_0(K\otimes_R\Lambda))\leq rank_{\Z}(K_0(B/I))$$
\end{Prop}

Proof. We shall first suppose that $\Lambda$ is indecomposable and that
$I\leq rad(B)$.

Let $S_1,S_2,\dots,S_n$ be representatives of the simple $B$-modules.
The projective cover $P_i$ of $S_i$ has then a composition series
where $soc(P_i)\simeq S_i$ and
$$rad(P_i)/soc(P_i)\simeq S_{i-1}\oplus S_{i+1}$$
for all $i\in \{2,3\dots,n-1\}$,
$$rad(P_1)/soc(P_1)\simeq S_{2}$$ and
$$rad(P_n)/soc(P_n)\simeq S_{n-1}.$$

Denote $\ol{B}:=B/I$. Since $I\leq rad(B)$,
we get $\ol{B}$ has the same number of simple modules, and
moreover, the simple
$\ol B$-modules and the simple $B$-modules coincide by the epimorphism
$B\lra\ol B$. Therefore
the projective indecomposable $\ol{B}$-modules are
$\ol{P}_i:=\ol{B}\otimes_BP_i$ for $i\in \{1,2,\dots,n\}$.
Moreover, $\ol{P}_i$ is the projective cover of $S_i$ as $\ol{B}$-module.

Extending $k$ if necessary, we may assume that
the field of fractions $K$ of $R$ is a splitting
field for $\Lambda$, since extending $K$ does not
decrease the rank of the Grothendieck group,
using the Noether-Deuring theorem.
Since $k$ is a splitting field for $B$ and for $B/I$,
and since $k\otimes_R\Lambda\simeq \ol{B}$, the Cartan matrix of $\ol{B}$
is symmetric (cf e.g. \cite[Proposition 4.2.11]{derbuch}).
Since the Cartan matrix of $B$ equals
$$C:=\left(\begin{array}{cccccc}
2&1&0&\dots&\dots&0\\1&2&1&\ddots&&\vdots\\
0&1&\ddots&\ddots&\ddots&\vdots\\
\vdots&\ddots&\ddots&\ddots&1&0\\
\vdots&&\ddots&1&2&1\\
0&\dots&\dots&0&1&2\end{array}\right),$$
we see that the composition length of $\ol P_i$ differs from
the composition length of $P_i$ by at most $1$. Indeed, if this would  not
be the case, then the composition matrix of $\ol B$ would be
decomposable into
at least two blocks and $\ol B$ would be decomposable as algebra.
But, since $\Lambda$ is indecomposable, so is $k\otimes_R\Lambda\simeq
\ol B.$ So, $I\leq soc(B)$.

Since $R$ is complete, we may assume that $B$ and
$\Lambda$ are both basic algebras.

Let $Q_i$ be the projective cover of $\ol P_i$ as $\Lambda$-module.
Hence, $k\otimes_RQ_i\simeq\ol P_i$. Since
$$dim_KHom_{K\Lambda}(KQ_i,KQ_j)=dim_RHom_{\Lambda}(Q_i,Q_j)=
dim_kHom_{\ol B}(\ol P_i, \ol P_j)$$
we know that $KQ_i$ and $KQ_j$ do not have a character in common if
$|i-j|>1$ and do have one character in common if $|i-j|=1$.
Since $dim_kHom_{\ol B}(\ol P_i, \ol P_i)\in\{1,2\}$,
the character of $KQ_i$ for all $i\in\{1,\dots,n\}$ is
a sum of at most two irreducible characters, and in case of
two characters these are non isomorphic. Now, since
$Hom_{K\Lambda}(KQ_i,KQ_j)=0$ if $|i-j|>1$, it follows that if
$$dim_KHom_{K\Lambda}(KQ_{i_0},KQ_{i_0})=1,$$ then $i_0\in\{1,n\}$.
Otherwise,
the character of $KQ_{i_0}$ would be a constituent of
$KQ_{i_0+1}$, of $KQ_{i_0}$ and of $KQ_{i_0-1}$, which implies
then $Hom_{K\Lambda}(KQ_{i_0+1},KQ_{i_0-1})\neq 0$. This would
give a contradiction.
This gives that $I$ equals either $S_1$ or $S_n$ or $S_1\oplus S_n$.

Suppose now that $S_1$ is a direct factor of $I$ (as left module)
and suppose $$rank_{\Z}(K_0(K\otimes_R\Lambda))> rank_{\Z}(K_0(B)).$$
We shall prove that $S_1$ is not a direct factor of $I$. By
symmetry, then neither $S_n$ is a direct factor of $I$, and therefore,
$I=0$.

Under these hypotheses, $dim(\ol P_1)=2$ and as a consequence also
$dim_R(Q_1)=2.$ So, for the Wedderburn components corresponding to
$K\otimes_RQ_1$ in $K\otimes_R\Lambda$ we have two possibilities.
Either $K\otimes_RQ_1$ is a sum of two one-dimensional characters
or $K\otimes_RQ_1$ is isomorphic to one two-dimensional character.
Since $K\otimes_R\Lambda$ admits at least $n+1$ irreducible characters,
$K\otimes_RQ_1$ must have two constituents. So, $K\otimes_RQ_1$
is a sum of two one-dimensional characters.

But now, let $\{e_i|i\in\{1,\dots,n\}\}$ be an orthogonal set
of primitive idempotents with $\Lambda e_i\simeq Q_i$. Then, since
$e_1$ and $e_2$ must be non zero on the common Wedderburn
component of $KQ_1$ and $KQ_2$, we get $e_2e_1\neq 0$.
This is a contradiction to the fact that $e_1$ and $e_2$ are orthogonal.

\medskip

We have to deal with the case $\Lambda$ being decomposable.
The structure of $B$ implies that in this case, $B/I$ is a direct
product of algebras we have dealt with in the earlier case,
and copies of $k$.
By induction on the number of simple modules of each indecomposable
factor the result holds for each of the pieces as well.
Summing up for all of these pieces, we get the desired result.

\medskip

Finally, we have to deal with the case that $I$ is not contained
in $rad(B)$. The same argument as for $\Lambda$ decomposable
applies here as well.

This proves the Proposition.
\phantom{x}\hfill\dickebox

\subsection{Proving that the
Baues-Dreckmann-Franjou-Pirashvili ring is an order}
\label{lifting}

Since $\hat\Z_p$ is a complete discrete valuation ring,
we may lift idempotents from $\Gamma_{\F_p}^p$ to $\Gamma_{\hat\Z_p}^p$.
Hence, there is an indecomposable direct factor $\Gamma_{\hat\Z_p,0}^p$ of
the rank one free module
$\Gamma_{\hat\Z_p}^p$ which maps surjectively to $\Gamma_{\F_p,0}^p$.
Let $T_0^p:=t(\Gamma_{\hat\Z_p,0}^p)$ be the torsion ideal in
$\Gamma_{\hat\Z_p,0}^p$ and define
$\Lambda_{\hat\Z_p,0}^p:=\Gamma_{\hat\Z_p,0}^p/T_0^p.$

\begin{Prop}\label{nplusone}
$\Lambda_{\hat\Z_p,0}^p$ is an order. Moreover,
$\hat\Q_p\otimes_{\hat\Z_p}\Lambda_{\hat\Z_p,0}^p$ is a direct
product of $p+1$ matrix rings over $\hat\Q_p$ and up to isomorphism
there are at most $p$ simple $\Lambda_{\hat\Z_p,0}^p$-modules.
\end{Prop}

Proof.
In fact, $\hat\Q_p\otimes_{\hat\Z_p}\Lambda_{\hat\Z_p}^p=
\hat\Q_p\otimes_{\hat\Z_p}\Gamma_{\hat\Z_p}^p$ and their
common module categories are equivalent to
the category of  polynomial functors $\hat Q_p-mod\lra\hat Q_p-mod$
of degree at most $p$ by Lemma~\ref{char0}.

By Friedlander-Suslin~\cite{friedlander}, the category of exact
degree $n$ polynomial
functors $\hat Q_p-mod\lra\hat Q_p-mod$ is equivalent to the category
of strict polynomial functors $\hat Q_p-mod\lra\hat Q_p-mod$ and this
category is
equivalent to the category of modules over the Schur algebra
$S_{\hat\Q_p}(n,n)$. Moreover, the category of strict polynomial
functors of degree at most $n$ is equivalent to the direct sum
of the category of strict polynomial functors of exact degree $m$ for
each $m\in\{0,1,\dots,n\}$. The Schur algebra $S_{\hat\Q_p}(p,p)$ is
split semisimple (cf Green~\cite{Green}) with exactly $p+1$
simple modules. This shows that
$\hat\Q_p\otimes_{\hat\Z_p}\Lambda_{\hat\Z_p,0}^p$ is a direct
product of $p+1$ full matrix rings over $\hat\Q_p$. Moreover, this shows
also that $\Lambda_{\hat\Z_p}^p$ is an order since it is by definition
torsion free and contains a basis of the semisimple algebra
$\hat\Q_p\otimes_{\hat\Z_p}\Lambda_{\hat\Z_p}^p$
(which is Morita equivalent
to $\prod_{i=0}^pS_{\hat\Q_p}(i,i)$).

In order to prove the second statement we just observe that
the number of simple objects in ${\mathcal A}_{\F_p}^p$ equals the number
of simple objects in ${\mathcal F}_{\F_p}^p$ by
Lemma~\ref{simplescorrespond}. Moreover, since $\Lambda_{\hat\Z_p,0}^p$
is a quotient of $\Gamma_{\hat\Z_p,0}^p$, every simple
$\Lambda_{\hat\Z_p,0}^p$-module induces a simple
$\Gamma_{\hat\Z_p,0}^p$-module. We know that $\Gamma_{\F_p,0}^p$
is a Brauer tree algebra with $p$ simple modules. Moreover,
$\F_p\otimes_{\hat\Z_p}\Gamma_{\hat\Z_p,0}^p\simeq \Gamma_{\F_p,0}^p$
and so, $\Gamma_{\hat\Z_p,0}^p$ admits $p$ simple modules. As a consequence
$\Lambda_{\hat\Z_p,0}^p$ admits at most $p$ simple modules.
This proves the proposition.\hfill\dickebox

\begin{Prop}\label{torsioniszero}
$t(\Gamma_{\F_p,0}^p)=0$ and therefore
$\Lambda_{\hat\Z_p,0}^p=\Gamma_{\hat\Z_p,0}^p$.
\end{Prop}

Proof. This is a consequence of Proposition~\ref{nplusone},
Theorem~\ref{gammapstructure} and Proposition~\ref{BBalg}.

Indeed, since $\F_p\otimes_{\hat\Z_p}-$ is right exact, the
epimorphism $$\Gamma_{\hat\Z_p,0}^p\lra \Lambda_{\hat\Z_p,0}^p$$
induces an epimorphism
$$B=\Gamma_{\F_p,0}^p\lra
\Lambda_{\hat\Z_p,0}^p\otimes_{\hat\Z_p}\F_p\;$$
with kernel $I$,
for $B$ being a Brauer tree algebra associated to a stem with $p$
edges and without exceptional vertex
(Theorem~\ref{gammapstructure}).
Since $t(\Gamma_{\hat\Z_p,0}^p)\subseteq rad(\Gamma_{\F_p,0}^p)$
by Proposition~\ref{nplusone}, $I\leq rad(B)$. Since
$B/I\simeq\Lambda_{\hat\Z_p,0}^p\otimes_{\hat\Z_p}\F_p\;$ for an
order $\Lambda_{\hat\Z_p,0}^p$, Proposition~\ref{BBalg} implies that in
this case $I=0$. Hence,
$\F_p\otimes_{\hat\Z_p}t(\Gamma_{\hat\Z_p,0}^p)=0$ and therefore,
$t(\Gamma_{\hat\Z_p,0}^p)=0$.
This proves the proposition.\hfill\dickebox

\subsection{Describing the order; the main result}
\label{determineorder}

We shall describe $\Lambda_{\hat\Z_p,0}^p$ and prove our main result.
For this purpose
we introduce some notation (cf \cite[Section 4.4]{derbuch}). Let

\noindent\unitlength1cm
\begin{center}\begin{picture}(6,.5)
\put(0,0){$\hat\Z_p$}
\put(.4,.15){\line(1,0){1}}
\put(1.5,0){$\hat\Z_p$}
\put(.9,.3){\scriptsize $p^i$}
\put(2.5,0){$:=\;\;\{(a,b)\in\hat\Z_p\times\hat\Z_p|\;a-b\in
p^i\hat\Z_p\}$}
\end{picture}\end{center}

and
\noindent\unitlength1cm
\begin{center}\begin{picture}(5,.3)
\put(0,0){$\hat\Z_p$}
\put(.4,.15){\line(1,0){1}}
\put(1.5,0){$\hat\Z_p$}
\put(2.5,0){$:=$}
\put(3.5,0){$\hat\Z_p$}
\put(3.9,.15){\line(1,0){1}}
\put(5,0){$\hat\Z_p\;\;$.}
\put(4.4,.3){\scriptsize $p$}
\end{picture}\end{center}

The following is the main result of our paper.

\begin{Theorem}\label{maintheorem}
Let ${\mathcal A}^p_{\hat\Z_p}$ be the category of at
polynomial functors from free abelian groups to $\hat\Z_p$-modules
and of degree at most $p$. Then,
${\mathcal A}^p_{\hat\Z_p}$ is equivalent to
$\Gamma_{\hat\Z_p}^p$-mod, where
$$\Gamma_{\hat\Z_p}^p:=(\prod_{1<n<p}\hat\Z_p)\times
(\prod_{\lambda\vdash p\mbox{ \scriptsize and} \;\lambda
\mbox{ \scriptsize not a hook}}\hat\Z_p)\times
\Lambda_{\hat\Z_p,0}^p$$
and where
\noindent
\begin{center}
\unitlength1cm\begin{picture}(12.5,3)
\put(0,2){$\Lambda_{\hat\Z_p,0}^p\simeq$}
\put(1.5,2){$\hat\Z_p\;\;\;\oplus$}
\put(2.8,2){$\left(\begin{array}{cc}\hat\Z_p&\hat\Z_p\\
(p)&\hat\Z_p\end{array}\right)\;\oplus$}
\put(5.8,2){$\left(\begin{array}{cc}\hat\Z_p&\hat\Z_p\\
(p)&\hat\Z_p\end{array}\right)\;\oplus$}
\put(8.5,2){$\dots\;\oplus$}
\put(9.6,2){$\left(\begin{array}{cc}\hat\Z_p&\hat\Z_p\\
(p)&\hat\Z_p\end{array}\right)\;\oplus$}
\put(12.3,2){$\hat\Z_p$}
\put(1.8,2.2){\line(6,1){1.4}}
\put(4.6,1.8){\line(3,1){1.4}}
\put(7.4,1.9){\line(2,1){.8}}
\put(9.1,1.9){\line(2,1){.9}}
\put(11.2,1.8){\line(4,1){1.1}}
\put(.1,.5){$=\left\{(d_0)\times\left(\prod_{j=1}^{p-1}
\left(\begin{array}{cc}a_j&b_j\\c_j&d_j\end{array}\right)\right)
\times (a_p)\;\vrule height3.1ex depth1.5ex\;
\forall j:
a_j,b_j,c_j,d_j\in\hat\Z_p;p|c_j;p|(d_j-a_{j-1})\right\}$}
\end{picture}
\end{center}
is a Green order with $p$ isomorphism classes of indecomposable
projective modules.
\end{Theorem}

\begin{Rem}
Roggenkamp described the orders $\Lambda$ which admit a set of lattices
with periodic projective resolutions encoded by a Brauer tree
(\cite{roggreen}, see also \cite{derbuch}).
Roggenkamp called these orders Green-orders and he described their
structure in great detail.
\end{Rem}

{\bf Proof of the theorem.}
The case $p\leq 3$ was done by Drozd. Hence we may suppose that $p\geq 5$.
Since $\F_p\otimes_{\hat\Z_p}
\Lambda_{\hat\Z_p,0}^p$ is a Brauer tree algebra,
there is a set of $\F_p\otimes_{\hat\Z_p}
\Lambda_{\hat\Z_p,0}^p$-modules having a periodic
projective resolution given by the Brauer tree of $\Lambda_{\F_p,0}^p $.
Lifting these projective resolutions to the order
$\Lambda_{\hat\Z_p,0}^p$ gives a periodic projective resolution
of certain $\Lambda_{\hat\Z_p,0}^p$-modules $M_i$.
These periodic resolutions are encoded by the same Brauer tree.
It remains to show that the modules $M_i$ are lattices. Actually, this is
automatic. Indeed, since the resolution is periodic, each module
$M_i$ is also a kernel of a differential, after a complete period of
the periodic projective resolution. Hence,
$\Lambda_{\hat\Z_p,0}^p$ is a Green order with Brauer tree
being a stem with $p$ edges and without exceptional vertex.

We have to show that the maximal overorder of the
Green order $\Lambda_{\hat\Z_p,0}^p$ is
a direct product of matrix rings over $\hat\Z_p$ and that the image
of $\Lambda_{\hat\Z_p,0}^p$ in each of the matrix rings is
a hereditary order.

The first part is clear since $\hat\Q_p$ is a splitting field
of $\Lambda_{\hat\Z_p,0}^p$, and $\Lambda_{\hat\Z_p,0}^p$ can be
embedded into a direct product of matrix rings over the ring of
integers in $\hat\Q_p$ (see e.g. \cite{MO}).
Let $e_1,e_2,\dots,e_{p+1}$ be a complete set of primitive
pairwise orthogonal idempotents of the center of
$\hat\Q_p\otimes_{\hat\Z_p}\Lambda_{\hat\Z_p,0}^p$.
Then,
$$\prod_{j=1}^{p-1}(\Lambda_{\hat\Z_p,0}^p\cdot e_i)\;\simeq\;
\hat\Z_p\times
\prod_{j=1}^{p-1}
\left(\begin{array}{cc}
\hat\Z_p&\hat\Z_p\\(p^{x_j})&\hat\Z_p
\end{array}\right)\times\hat\Z_p
$$
for some $x_j\in{\mathbb N}\setminus\{0\}$.
Moreover, since $\Lambda_{\hat\Z_p,0}^p\otimes_{\hat\Z_p}\F_p$
is a Brauer tree algebra without exceptional vertex,
$x_1=x_2=\dots=x_{p-1}$ and as a consequence,
if one of the matrix rings is hereditary,
all of them are hereditary. The structure theory of Green orders
(cf Roggenkamp \cite{roggreen}; see also \cite[Section 4.4]{derbuch})
and of hereditary orders,
(cf e.g. Reiner \cite{MO}) then gives the statement.

Define a functor
$$Hom_{{\mathcal A}_{\hat\Z_p}^p}\left(\hat\Z_p
\otimes_\Z(\bigotimes_{j=1}^pid),-\right):
{\mathcal A}_{\hat\Z_p}^p\lra \hat\Z_p\Sn_p-mod$$
where we use again that the functor
$\hat\Z_p\otimes_\Z(\bigotimes_{j=1}^pid):\Z-free\lra \hat\Z_p-mod$
carries a natural $\hat\Z_p$-linear $\Sn_p$ action. Denote for notational simplicity
$E:=Hom_{{\mathcal A}_{\hat\Z_p}^p}
\left(\hat\Z_p\otimes_\Z(\bigotimes_{j=1}^pid),-\right)$.
By Lemma~\ref{macdonaldlinearisation} this functor is just
the $p-1$-th cross effect, and by Lemma~\ref{degreeunderext}
this functor $E$ is exact.

Since $E$ is exact, and since by definition $E$ is represented by
$P:=\hat\Z_p\otimes_\Z(\bigotimes_{j=1}^pid)$, this object $P$ is
projective. Let $e$ be an idempotent in $\Gamma_{\hat\Z_p}^p$
which corresponds to the projective indecomposable
$\Gamma_{\hat\Z_p}^p$-modules which occur in $P$.
Then, replacing ${\mathcal A}_{\hat\Z_p}^p$ by $\Gamma_{\hat\Z_p}^p-mod$
the functor $Hom_{{\mathcal A}_{\hat\Z_p}^p}\left(\hat\Z_p
\otimes_\Z(\bigotimes_{j=1}^pid),-\right)$ becomes
the functor $E:\Gamma_{\hat\Z_p}^p-mod\lra \hat\Z_p\Sn_p-mod$
and $E$ is just multiplication by $e$.

We need to show that
$End_{{\mathcal A}_{\hat\Z_p}}(P)\simeq \hat\Z_p\Sn_p$, where the
action is given by permutation of components in the tensor product.
Once this is done, we know that
$End_{{\mathcal A}_{\hat\Z_p}}(P)\simeq
e\cdot\Gamma_{\hat\Z_p}^p\cdot e\simeq \hat\Z_p\Sn_p$
and we observe that for all idempotents in $\Gamma_{\hat\Z_p}^p$
we get that this product $e\cdot\Gamma_{\hat\Z_p}^p\cdot e$ is again a
product of Green orders with the same order of congruences.
Since $\hat\Z_p\Sn_p$ is a Green order with congruences modulo $p$ only,
we get that $x=1$.

\begin{Claim}
$End_{{\mathcal A}_{\hat\Z_p}}
(\hat\Z_p\otimes_\Z(\bigotimes_{j=1}^pid))\simeq \hat\Z_p\Sn_p$
\end{Claim}

Proof.
The proof given by Piriou-Schwartz \cite[Lemma 1.9]{PS}
of the corresponding statement for ${\mathcal F}_{\F_p}$
carries over literally. For the reader's convenience we recall the
(short) arguments.

Given an $x=\sum_{\sigma\in\Sn_p}x_\sigma\sigma\in\hat\Z_p\Sn_p$, then
associate to this $x$ the natural transformation $\eta_x$ in
$End_{{\mathcal A}_{\hat\Z_p}}
\left(\bigotimes_{j=1}^p\left(\hat\Z_p\otimes_\Z id\right)\right)$ given by
$$v_1\otimes v_2\otimes\dots\otimes v_n\mapsto
\sum_{\sigma\in\Sn_p}x_\sigma (v_{\sigma^{-1}(1)}
\otimes\dots\otimes v_{\sigma^{-1}(n)}).$$
Inversely, a natural transformation $\eta$ in
$End_{{\mathcal A}_{\hat\Z_p}}
\left(\bigotimes_{j=1}^p\left(\hat\Z_p\otimes_\Z id\right)\right)$
is determined by its
value on $\Z^n$. Fix a basis $\{e_1,\dots, e_n\}$ of $\Z^n$.
The image of $e_1\otimes e_2\otimes\dots\otimes e_n$ under $\eta_{\Z^n}$
can be uniquely
written as $x_\eta\cdot (e_1\otimes e_2\otimes\dots\otimes e_n)$ for an
$x_\eta\in\hat\Z_p\Sn_n$. The two mappings $x\mapsto \eta_x$ and
$\eta\mapsto x_\eta$ are mutually inverse and
obviously ring homomorphisms. \hfill\dickebox

This proves the theorem.
\phantom{x}\hfill\dickebox

\begin{Rem}
The
Schur algebra $S_{\hat\Z_p}(p,p)$ is a classical order which was
completely described by K\"onig in \cite{cyclotomic}.

\noindent
\begin{center}
\unitlength1cm\begin{picture}(12.5,1.1)
\put(0,.5){$S_{\hat\Z_p}(p,p)'\simeq$}
\put(2.5,.5){$\left(\begin{array}{cc}\hat\Z_p&\hat\Z_p\\
(p)&\hat\Z_p\end{array}\right)\;\oplus$}
\put(5.5,.5){$\left(\begin{array}{cc}\hat\Z_p&\hat\Z_p\\
(p)&\hat\Z_p\end{array}\right)\;\oplus$}
\put(8.5,.5){$\dots\;\oplus$}
\put(9.5,.5){$\left(\begin{array}{cc}\hat\Z_p&\hat\Z_p\\
(p)&\hat\Z_p\end{array}\right)\;\oplus$}
\put(12.5,.5){$\hat\Z_p$}
\put(4.3,.3){\line(3,1){1.7}}
\put(7.3,.5){\line(2,1){.8}}
\put(9.1,.3){\line(2,1){.9}}
\put(11.3,.3){\line(4,1){1.1}}
\end{picture}
\end{center}
where we denote by $S_{\hat\Z_p}(p,p)'$ the basic algebra of the Schur
algebra $S_{\hat\Z_p}(p,p)$.

Now, any strict polynomial functor induces a polynomial functor.
So, composing further to the Green order lifting the principal block of
the group ring of the symmetric group, we get an induced functor
$$
S_{\hat\Z_p}(p,p)'-mod\lra {\mathcal A}^p_{\hat\Z_p}\stackrel{\phi}{\lra}
e\cdot \Lambda\cdot e-mod.
$$
Since the functor $\phi$ is induced by the Schur functor,
this composed map is induced by the natural
embedding of $e\cdot \Lambda\cdot e\hookrightarrow
S_{\hat\Z_p}(p,p)'$.
\end{Rem}

\section{Identifying the lattices as functors}

We shall identify the indecomposable functors of ${\mathcal A}_{\hat
\Z_p}^p$ which correspond to indecomposable
$\Gamma_{\hat\Z_p}$-lattices. We call such polynomial functors
'polynomial lattices'.

The structure of $\Gamma_{\hat\Z_p,0}$ implies that there are
exactly $3p-2$ indecomposable $\Gamma_{\hat\Z_p,0}$ lattices.
Indeed, the indecomposable lattices are
the $p$ projective indecomposable modules $P_1, P_2, \dots, P_p$,
the $p-1$ kernels of any fixed non zero homomorphism $P_i\lra
P_{i+1}$ for $i=1,2,\dots,p-1$, as well as the $p-1$
kernels of any fixed non
zero homomorphism $P_i\lra P_{i-1}$ for $i=2,3,\dots,p$.
Therefore, there are exactly $3p-2$ indecomposable 'lattices' of
exactly degree $p$ polynomial functors in ${\mathcal A}_{\hat\Z_p}^p$.

Moreover, the proof of Theorem~\ref{maintheorem} shows that
exactly the projective indecomposable $\Gamma_{\hat\Z_p}$-modules
of degree $d\in\{2,\dots,p-1\}$ will give rise to indecomposable
lattices. Denote by $\rho(k)$ the number of partitions of $k$ into
non zero integers, we get the following corollary to
Theorem~\ref{maintheorem}.

\begin{Cor}
Up to isomorphism there are exactly
$3p-2+\sum_{k=2}^{p-1}\rho(k)$ indecomposable
polynomial lattices in ${\mathcal A}_{\hat\Z_p}^p$ and
$p+\sum_{k=2}^{p-1}\rho(k)$ of them are projective, while $2(p-1)$ of
them are not projective. The non projective polynomial lattices
are kernels of mappings between projective indecomposable
polynomial functors.
\end{Cor}


\begin{thebibliography}{88}

\bibitem{ARS}
Maurice Auslander, Idun Reiten and Sverre Smal\o{}, {\sc Represenation
Theory of Artin Algebras}, Cambridge University Press, Cambridge 1995

\bibitem{Baues}
Hans-Joachim Baues, {\em Quadratic functors and metastable
homotopy}, Journal of pure and applied algebra {\bf 91} (1994)
49-107.

\bibitem{BDFP}
Hans-Joachim Baues, Winfried Dreckmann, Vincent Franjou and Teimuraz
Pirashvili, {\em Foncteurs
polynomiaux et foncteurs de Mackey non lin\'eaires}, Bulletin de la
Soci\'et\'e de Math\'ematiques de France {\bf 129} (2001) 237-257.

\bibitem{Benson}
David Benson, {\sc Representations and Cohomology} Vol 1. Cambridge University Press, Cambridge 1991.

\bibitem{Djament}
Aur\'elien Djament, {\em Sur l'homologie des groupes unitaires à coefficients polynomiaux. } Journal of K-Theory {\bf 10} (2012) 87-139.

\bibitem{DjamentVespa}
Aur\'elien Djament and Christine Vespa, {\em Sur l'homologie des groupes orthogonaux et symplectiques à coefficients tordus}, Annales Scientifiques
de l'\'Ecole Normale Sup\'erieure {\bf 43} (2010) 395-459.

\bibitem{DjamentVespa2}
Aur\'elien Djament and Christine Vespa, {\em Sur l'homologie des groupes d'automorphismes des groupes libres \`a coefficients polynomiaux}, preprint 2013,
arxiv: 1210.4030v2

\bibitem{drozdquadratic}
Yuri Drozd, {\em Finitely generated quadratic modules},
manuscripta mathematica {\bf 104} (2001) 239-256.

\bibitem{drozdcubic}
Yuri Drozd, {\em On cubic functors}, Communications in Algebra
{\bf 31} (3) (2003) 1147-1173.

\bibitem{eilenberg}
Samuel Eilenberg and Saunders Maclane, {\it On the groups
$H(\Pi,n)$, II}; Annals of Mathematics {\bf 70} (1954) 49-139.

\bibitem{Franjouext}
Vincent Franjou, {\em Extensions entre puissances ext\'erieures et entre puissances sym\'etriques}, Journal of Algebra {\bf 179} (1996) 501-522.

\bibitem{FLS}
Vincent Franjou, Jean Lannes and Lionel Schwartz,
{\em Autour de la cohomologie de MacLane des corps finis},
Inverntiones Mathematicae {\bf 115} (1994) 513-538.

\bibitem{FP}
Vincent Franjou and Teimuraz Pirashvili,
{\em On MacLane Cohomology for the ring of integers}, Topology
{\bf 37} (1998) 109-114

\bibitem{skt}
Vincent Franjou and Teimuraz Pirashvili,
{\em Stable $K$-theory is bifunctor homoogy (after A. Scorichenko)},  in
{\sc Franjou, Friedlander, Pirashviliand Schwartz:
Rational representations, the Steenrod algebra and functor homology},
Panoramas et Synth\`eses 16, Soci\'et\'e Math\'ematique de France, 2003.

\bibitem{FFSS}
Vincent Franjou, Eric M. Friedlander, Andrei Scorichenko and Andrei Suslin,
{\em General linear and functor cohomology over finite fields},
Annals of Mathematics {\bf 155} (1999) 663-728.

\bibitem{friedlander}
Eric M. Friedlander and Andrei Suslin, {\em Cohomology of finite
group schemes over a field}, Inventiones mathematicae {\bf 127} (1997)
209-270.

\bibitem{Green}
J. A. Green, {\sc  Polynomial representations of $GL_n$}, Springer
Lecture Notes in Mathematics 830 (1980).

\bibitem{HPV}
Manfred Hartl, Teimuraz Pirashvili and Christine Vespa, {\em Polynomial Functors from Algebras over a Set-Operad and Non-Linear Mackey Functors},
preprint (2012) arxiv:1209.1607v2

\bibitem{HLS}
Hans-Werner Henn, Jean Lannes and Lionel Schwartz, {\em The category of
unstable modules and unstable algebras over the Steenrod algebra modulo
nilpotent objects}, American Journal of Mathematics {\bf 115}
(1993) 1053-1106.

\bibitem{cyclotomic}
Steffen K\"onig, {\em Cyclotomic Schur algebras and blocks of cyclic defect},
Canadian Mathematical Bulletin {\bf 43} (2000) 79-86.

\bibitem{derbuch}
Steffen K\"onig and Alexander Zimmermann, {\sc Derived equivalences
for group rings}, Springer Lecture Notes in Mathematics 1685 (1998).

\bibitem{Kuhn1}
Nicholas Kuhn, {\em Generic Representations of the finite general linear
groups and the Steenrod algebra: I}, American Journal of Mathematics
{\bf 116} (1993) 327-360.

\bibitem{KuhnII}
Nicholas Kuhn, {\em Generic Representations of the finite general linear
groups and the Steenrod algebra: II}, K-theory {\bf 8} (1994)
395-428.


\bibitem{kuhnsurvey}
Nicholas J. Kuhn, {\em The generic representation theory of finite
fields: A survey of basic structure}, in: Infinite length modules,
(Bielefeld 1998) 193-212, Trends in Mathematics, Birkh\"auser,
Basel (2000).

\bibitem{stratifying}
Nicholas Kuhn, {\em A stratification of generic representation theory and
generalized Schur algebras}, $K$-Theory {\bf 26} (2002) 15-49.

\bibitem{Maclane}
Saunders Maclane, {\sc Categories for the working mathematician}, second edition,
Springer Verlag, Heidelberg 1997

\bibitem{Pirashvilirussian}
Teimuraz Pirashvili, {\em Polynomial functors} (Russian, English summary). Trudy Tibliss. Mat. Inst. Razmadze Akad. Nauk Gruzin. SSR {\bf 91} (1988) 5566.

\bibitem{pirashvilibourbaki}
Teimuraz Pirashvili, {\em Polynomial functors over finite fields},
S\'eminaire Bourbaki, Volume 1999/2000, expos\'e no 865-879;
Ast\'erisque {\bf 276} (2002) 369-388.

\bibitem{pirashvililecturenotes}
Teimuraz Pirashvili, {\em Introduction to functor homology}, in
{\sc Franjou, Friedlander, Pirashviliand Schwartz:
Rational representations, the Steenrod algebra and functor homology},
Panoramas et Synth\`eses 16, Soci\'et\'e Math\'ematique de France, 2003.

\bibitem{piriouthese}
Laurent Piriou, {\em Extensions entre foncteurs de la cat\'egorie
des espaces vectoriels sur le corps premier \`a $p$ \'el\'ements
dans elle-m\^eme.}, th\`ese de doctorat universit\'e Paris 7
(1995).

\bibitem{PS}
Laurent Piriou and Lionel Schwartz, {\em Extension de foncteurs
simples}, $K$-Theory {\bf 15} (1998) 269-291.

\bibitem{Chrysostomos}
Chrysostomos Psaroudakis and Jorge Vot\'oria, {\em Recollements of module categories}, arXiv:1304.2692v1.

\bibitem{MO}
Irving Reiner, {\sc Maximal Orders}, Academic Press 1975.

\bibitem{roggreen}
Klaus W. Roggenkamp, {\em Blocks of cyclic defect and Green--orders,}
Communications in Algebra {\bf 20} (1992) 1715--1734.

\bibitem{Touze}
Antoine Touz\'e, {\em Cohomologie Rationnelle du Groupe Lin\'eaire et
Extensions de Bifoncteurs}, Th\`ese de doctorat de
l'universit\'e de Nantes (2008).

\bibitem{TouzevanderKallen}
Antoine Touz\'e and Wilberd van der Kallen, {\em Bifunctor cohomology and Cohomological finite generation for reductive groups},  Duke Mathematical Journal {\bf 151}  (2010),  251-278.

\bibitem{Vespa0}
Christine Vespa, {\em Generic representations of orthogonal groups: the mixed functors.} Algebraic and Geometric Topology {\bf 7} (2007) 379-410.

\bibitem{Vespa1}
Christine Vespa, {\em Generic representations of orthogonal groups: the functor category ${\mathfrak F}_{\rm quad}$.} Journal of Pure and Applied Algebra {\bf 212} (2008) 1472-1499.

\bibitem{Vespa2}
Christine Vespa, {\em Generic representations of orthogonal groups: projective functors in the category ${\mathfrak F}_{\rm quad}$.} Fund. Math. {\bf 200} (2008) 243-278.

\bibitem{unstable}
Lionel Schwartz, {\sc Unstable modules over the Steenrod algebra and
Sullivan's fixed point conjecture}, Chicago Lectures in Mathematics,
University of Chicago Press 1994.




\end{thebibliography}
\end{document}